\documentclass[12pt]{article}

\usepackage{arxiv}

\usepackage[utf8]{inputenc} 
\usepackage[T1]{fontenc}    
\usepackage{hyperref}       
\usepackage{url}            
\usepackage{booktabs}       
\usepackage{amsfonts}       
\usepackage{nicefrac}       
\usepackage{microtype}      
\usepackage{graphicx}
\usepackage{doi}
\usepackage{booktabs,graphicx,makecell,siunitx,eqparbox}
\usepackage{amsthm}
\usepackage{amsmath}
\usepackage{algorithm}
\usepackage{amssymb}
\usepackage{latexsym}
\usepackage{subcaption}
\usepackage{indentfirst}
\usepackage{caption}
\usepackage[capitalise]{cleveref}       
\usepackage{multirow}
\usepackage{diagbox}
\usepackage{todonotes}
\numberwithin{equation}{section}

\title{Numerical study on hyper parameter settings for neural network approximation
	to partial differential equations}

\newif\ifuniqueAffiliation
\uniqueAffiliationtrue

\ifuniqueAffiliation 

\author{
	Hee Jun Yang 
	\thanks{
		National Institute for Mathematical Sciences, Korea. {Email:yangheejun1009@nims.re.kr.}
		The research of Hee Jun Yang is supported by National Institute for
		Mathematical Sciences (NIMS) grant funded by the Korea
		government (MSIT) (No. B25810000)
	}
	\And
	Alexander Heinlein 
	\thanks{
		Delft Institute of Applied Mathematics, Delft University of Technology, The Netherlands.
		{Email:a.heinlein@tudelft.nl.}
	}
	\And
	Hyea Hyun Kim 
	\thanks{Department of Applied Mathematics and Institute of Natural Sciences, Kyung Hee University,
		Korea. {Email:hhkim@khu.ac.kr.}
		The research of Hyea Hyun Kim is
		supported by the National Research Foundation of Korea(NRF) grants
		funded by NRF-2022R1A2C100388511.}
}

\fi

\renewcommand{\shorttitle}{}

\hypersetup{
	pdftitle={A template for the arxiv style},
	pdfsubject={q-bio.NC, q-bio.QM},
	pdfauthor={David S.~Hippocampus, Elias D.~Striatum},
	pdfkeywords={First keyword, Second keyword, More},
}

\begin{document}
\maketitle

\begin{abstract}
	Approximate solutions of partial differential equations (PDEs) obtained by neural networks are highly affected by hyper parameter settings. For instance, the model training strongly depends on loss function design, including the choice of weight factors for different terms in the loss function, and the sampling set related to numerical integration; other hyper parameters, like the network architecture and the optimizer settings, also impact the model performance. On the other hand, suitable hyper parameter settings are known to be different for different model problems and currently no universal rule for the choice of hyper parameters is known. 
	
	In this paper, for second order elliptic model problems, various hyper parameter settings are tested numerically to provide a practical guide for efficient and accurate neural network approximation. While a full study of all possible hyper parameter settings is not possible, we focus on studying the formulation of the PDE loss as well as the incorporation of the boundary conditions, the choice of collocation points associated with numerical integration schemes, and various approaches for dealing with loss imbalances will be extensively studied on various model problems; in addition to various Poisson model problems, also a nonlinear and an eigenvalue problem are considered.
\end{abstract}


\section{Introduction}

Recent advances in neural networks (NNs) have led to growing research efforts into their application in engineering and scientific applications. A particularly popular approach involves using NNs to discretize partial differential equations (PDEs), offering an alternative to classical numerical methods such as finite differences, finite elements, and finite volumes. First variants of methods where PDEs are incorporated into neural network training via the loss function were already introduced in seminal works from the 1990s~\cite{dissanayake_neural-network-based_1994,lagaris_artificial_1998}, shortly after key mathematical breakthroughs in the theory of neural networks, including the establishment of their universal approximation properties~\cite{cybenko_approximation_1989}. While many modern approaches have been developed during the past few years, physics-informed neural networks (PINNs)~\cite{raissi2019} and the deep Ritz method~\cite{yu2018deep} have been particularly successful. The whole class of methods is often generally referred to as \textit{physics-informed}; cf.~\cite{sirignano2018,kharazmi_hp-vpinns_2021,nguyen-thanh_deep_2020} for other related approaches. We also refer to the review articles~\cite{blechschmidt_three_2021,cai_physics-informed_2021,karniadakis_physics-informed_2021,cuomo_scientific_2022,toscano_pinns_2024,raissi_physics-informed_2024} for a more complete literature overview.

Physics-informed neural network approaches are generally easy to implement using state-of-the-art machine learning frameworks with automatic differentiation support, for instance, Tensorflow~\cite{abadi_tensorflow_2016}, PyTorch~\cite{paszke_pytorch_2019}, and Jax~\cite{bradbury_jax_2018}, without explicitly requiring a computational mesh~\cite{berg2018}. Moreover, they show great potential for addressing challenges such as incorporating observational data~\cite{kissas2020} or high-dimensional, inverse, and uncertainty quantification problems~\cite{weinan2017,Yang2019,yang2021b}. However, they also exhibit certain weaknesses that hinder their success in practical applications. In particular, they are difficult to train, and standard neural network optimizers are far from competitive with optimized numerical solvers used in classical numerical discretizations for most types of forward problems. This challenge appears closely related to the spectral bias or frequency principle of neural networks~\cite{rahaman_spectral_2019,xu_overview_2022}, i.e., the tendency of neural networks to approximate low-frequency components of functions more easily than high-frequency components. One possible explanation involves the spectral decomposition of the neural tangent kernel (NTK)\cite{jacot_neural_2018}, which provides insights into the convergence behavior of neural network training; see, for example,~\cite{NTK-theory} for a discussion in the context of PINNs. The spectral bias also makes the neural network training particularly difficult for multiscale and multifrequency problems. Another perspective on the failure of the training of PINN models is given in~\cite{basir2023investigating}. Successful approaches to improve the performance of PINNs involve adaptive weighting~\cite{mcclenny2023self} and sampling methods~\cite{lu_deepxde_2021}, advanced optimization techniques~\cite{muller_achieving_2023}, multi-stage~\cite{wang_multi-stage_2024} or multifidelity training approaches~\cite{howard_stacked_2024}, or domain decomposition-based approaches~\cite{D3M,dolean2024multilevel,yang_iterative_2024,jang2024partitioned,sun2024domain}.

Another major drawback of NN-based discretizations for PDEs is that the training and approximation properties strongly depend on the hyper parameter settings, including but not limited to the network architecture,
the loss function, the sampling of the training points, and the optimizer employed for training. Moreover, it is often observed that the optimal choice of parameters is highly problem-dependent. A study detailing some state-of-the-art choices in 2023 can be found in~\cite{wang2023expert}. Similarly, the model performance may strongly depend on the initialization of the trainable network parameters. These strong sensitivities often make it extremely difficult to reproduce results, once the problems settings are even varied only slightly. Nonetheless, many previous works did not investigate the sensitivity of the methods with respect to hyper parameter choices and network initialization.

In this paper, we present a detailed study of the performance of the two most popular physics-informed neural network approaches for approximating the solutions of PDEs, that is, PINNs and the deep Ritz method, depending on the initialization of neural network parameters and various hyper parameter choices. In particular, we will consider:
\begin{itemize}
	\item \textbf{PDE loss term formulations}: PINNs and deep Ritz method
	\item \textbf{Sampling schemes}: different from~\cite{wu2023comprehensive}, which compares different non-adaptive and residual-based Monte-Carlo sampling strategies, we focus on a comparison with Gaussian numerical integration schemes
	\item \textbf{Schemes for balancing the PDE and boundary loss terms}, including: constant and self-adaptive weights~\cite{mcclenny2023self} and an augmented Lagrangian approach~\cite{Son2023}
	\item \textbf{Neural network structure}: varying activation functions, Ansatz for hard enforcement of boundary conditions, and Fourier feature embedding~\cite{tancik_fourier_2020}
	\item \textbf{Optimizers}: Adam (adaptive moments)~\cite{adam} and LBFGS (limited-memory Broyden--Fletcher--Goldfarb--Shanno)~\cite{liu1989limited} algorithms
\end{itemize}

Our goal is to supplement the study of~\cite{wang2023expert} and come up with guidelines for the hyper parameter settings for neural network-based discretization methods depending on the model problem complexity. 
Our work is not a repetition of~\cite{wang2023expert} but extends its scope from only PINNs to also include the deep Ritz method and considering additional techniques; notably, for some challenging examples, we indeed observe advantages of the deep Ritz method in terms of the approximate solution accuracy and the training time.

This paper is organized as follows. In Section~\ref{model:method}, we introduce the model problems as well as the PINN and deep Ritz methods that form the basis of our numerical studies. Furthermore, we introduce some of the approaches to be compared, including sampling schemes based on Monte--Carlo and Gaussian numerical integration as well as different formulations for treating boundary conditions. Then, we introduce the detailed settings of our numerical experiments and list all employed hyper parameters in Section~\ref{numerics:settings}. In Section~\ref{numerics:result}, we report the results of our numerical experiments in order to come up with guidelines for the hyper parameter settings, depending on the complexity of the considered model problems. Then, in Section~\ref{sec:DeepRitz:numerics}, we present results for some more challenging three-dimensional, nonlinear, and eigenvalue problems. Finally, we add some further remarks and draw conclusions in Section~\ref{sec:Conclude}.

\section{Model problems and neural network approximation}\label{model:method}

In this section, we introduce the two-dimensional Poisson model problems that we will consider for the main part of our numerical studies; additional three-dimensional, nonlinear, and eigenvalue model problems will be introduced and studied in~\cref{sec:DeepRitz:numerics}. Afterwards, we will also introduce the neural network approximation schemes along with the hyper parameters investigated for their impact on the solution accuracy and efficiency.

\subsection{Poisson model problems} \label{sec:model_problems}

We consider the following Poisson problem on a unit square domain $\Omega=(0 \, 1)^2$,
\begin{equation}\label{model:poisson}
	\begin{split}
		-\nabla \cdot ( \nabla u ) &= f\quad \text{ in } \Omega,\\
		u&= g \quad \text{ on } \partial \Omega,
	\end{split}
\end{equation}
where we assume that the solution $u$ exists uniquely for the given functions $f$ and $g$.

In order to study various hyper parameter settings for the neural network models, we will consider the following variations of~\cref{model:poisson} and which are characterized by exact solutions with different complexity.

\paragraph{Example 1} Smooth and oscillatory solution with a positive integer $k$:
\begin{equation}\label{ex1}
	u(x,y)=\sin(k \pi  x)\sin( k \pi y).
\end{equation}

\paragraph{Example 2} Multi-frequency component solution with a positive integer $N$:
\begin{equation}\label{ex2}
	u(x,y)=\frac{1}{N}\sum_{\ell=1}^N \sin(2^\ell \pi x)\sin(2^\ell \pi y).
\end{equation}

\paragraph{Example 3} High contrast and oscillatory interior layer solution with $A>0$ and $\varepsilon>0$:
\begin{equation}\label{ex3}
	u(x,y)=Ax(1-x)y(1-y) \sin\left(\frac{(x-0.5)(y-0.5)}{\varepsilon}\right),
\end{equation}
where a large value $A$ and a small value $\varepsilon$ are considered.

\begin{figure}[ht!]
	\begin{center}
		\includegraphics[width=0.32\textwidth]{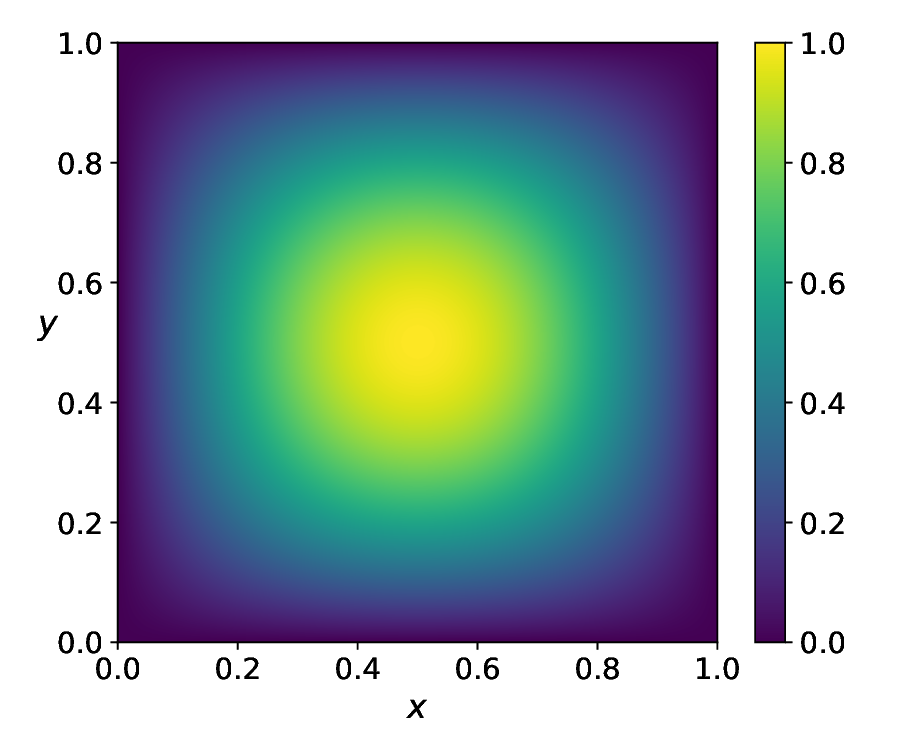}\;
		\includegraphics[width=0.32\textwidth]{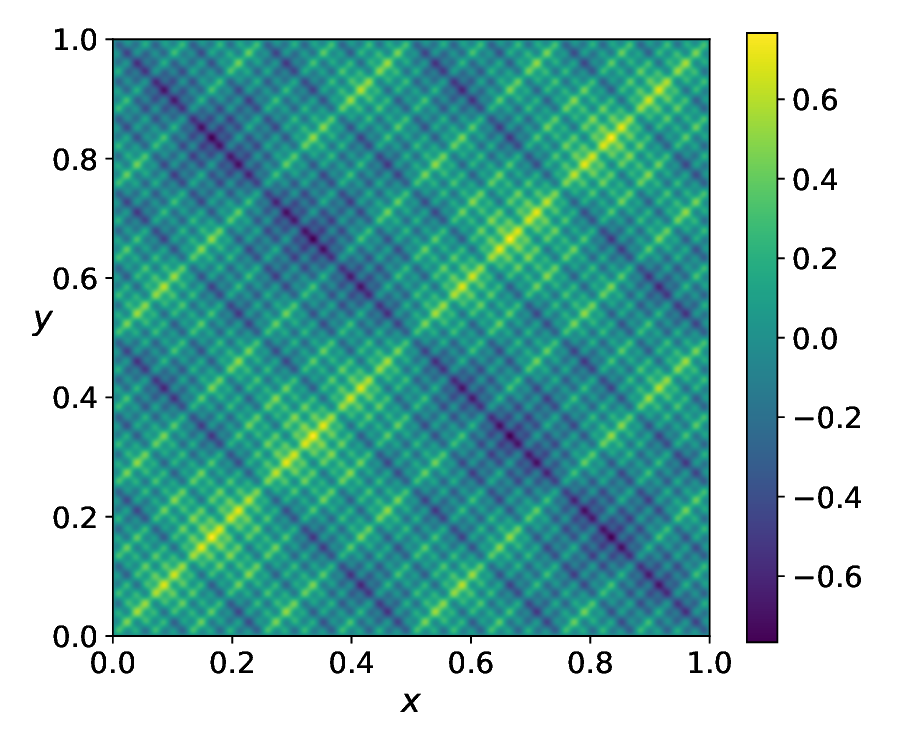}\;
		\includegraphics[width=0.32\textwidth]{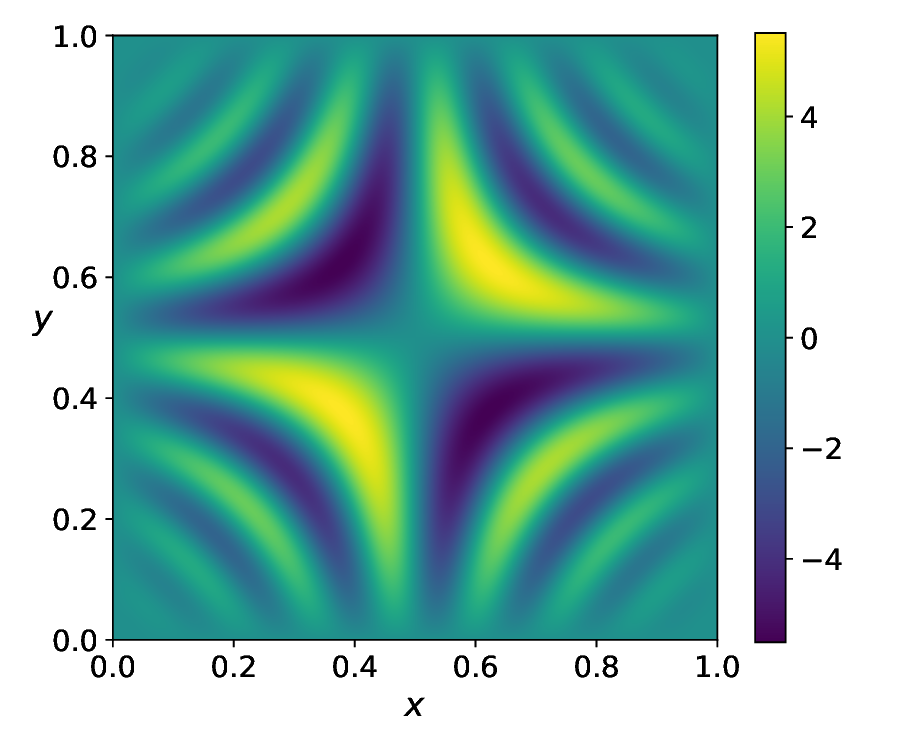}
	\end{center}
	\vskip-.7truecm
	\caption{Examples~\eqref{ex1}--\eqref{ex3}: solution plots for $k=1$ (left), $N=6$ (middle), and $A=100$ ans $\varepsilon=0.01$ (right), respectively.}
	\label{fig:examples}
\end{figure}

In Figure~\ref{fig:examples}, exemplary plots of the solutions of the three test examples with values $k=1$, $N=6$, and $A=100$ and $\varepsilon=0.01$, respectively, are presented.

\subsection{Neural network approximation} \label{sec:nn_approximation}

In order to approximate the solution of the model problems, we employ a neural network function $U({\bf x};\theta)$, where $\theta$ denotes the parameters of the neural network function. In training the neural network, the parameters are determined so as to satisfy the given differential equation and boundary condition of the model problem. In particular, a loss function related to the model problem is formed and the parameters are trained to minimize the loss function. In physics-informed neural networks (PINNs)~\cite{raissi2019}, the following form of loss function is introduced
\begin{equation} \label{loss:PINN}
	J_{P,M} (\theta):=\frac{w_I}{|X(\Omega)|} \sum_{ {\bf x} \in X(\Omega)} (\nabla \cdot \nabla U({\bf x};\theta)+f({\bf x}))^2 + \frac{w_B}{|X(\partial \Omega)|} \sum_{ {\bf x} \in X(\partial \Omega)} (U({\bf x};\theta)-g({\bf x}))^2,
\end{equation}
and the parameters $\theta$ in the neural network function $U({\bf x};\theta)$ are trained to minimize the loss value $J_{P,M}(\theta)$ in order to satisfy the differential equation and the boundary (or initial) condition of the model problem.
In the above, $X(A)$ denotes a set of training sampling points chosen from the domain $A$, $|X(A)|$ denotes the number of points in the set $X(A)$, and $w_I$ and $w_B$ are weight factors for the corresponding loss terms.
These weight factors are introduced to deal with the imbalance in the different terms in the loss function; see~\cite{NTK-theory}. The choices of the sampling sets $X(\Omega)$ and $X(\partial \Omega)$ as well as of the weight factors $w_I$ and $w_B$ are important hyper parameters of the PINN algorithm as they strongly impact the performance of the trained neural network model.

We note that the loss function $J_{P,M}(\theta)$ in~\cref{loss:PINN} is obtained from the Monte--Carlo approximation to the integrals of the residual of the differential equation and the boundary error,
$$
	J_P(\theta)
	:=
	w_I \int_{\Omega} (\nabla\cdot\nabla U({\bf x};\theta)+f({\bf x}))^2 \, d{\bf x} + w_B \int_{\partial \Omega} (U({\bf x};\theta)-g({\bf x}))^2 \, ds({\bf x}).
$$
We employ the subscript $M$ in the notation $J_{P,M}(\theta)$ to indicate that Monte--Carlo integration is used to approximate the integrals in the loss function $J_P(\theta)$, and the subscript $P$ in $J_P(\theta)$ to stress that the PDE loss function is formed by the PINN approach.

Instead of including the boundary condition in the loss function, {\it ansatz functions} $A({\bf x})$  and $G({\bf x})$ can be used to form the neural network function $\widetilde{U}({\bf x};\theta)$ so as to enforce the boundary (or the initial) condition strongly, i.e., as hard constraints,
\begin{equation} \label{eq:hard_constraints}
	\widetilde{U}(x;\theta):=A({\bf x})+G({\bf x})U({\bf x};\theta),
\end{equation}
where $A({\bf x})$ satisfies the boundary condition, $A({\bf x})=g({\bf x})$ on $\partial \Omega$ and $G({\bf x})=0$ on $\partial \Omega$; cf.~\cite{lagaris_artificial_1998}. The parameter $\theta$ in $\widetilde{U}(x;\theta)$ can then be trained to minimize
the loss function with only the differential equation term,
$$
	J_{P_I,M} (\theta)
	:=
	\frac{1}{|X(\Omega)|} \sum_{ {\bf x} \in X(\Omega)} (\nabla \cdot \nabla \widetilde{U}({\bf x};\theta)+f({\bf x}))^2,
$$
without the need for dealing with the weight factors to different terms in the loss function. Here, the subscript $P_I$ in $J_{P_I,M}(\theta)$ is employed to indicate the use of hard boundary constraints, using only the differential equation term defined on the domain interior in the PINN formulation. In contrast to hard boundary constraints, the penalty formulation of the boundary conditions in~\cref{loss:PINN} is also referred to as soft boundary constraints.

Other successful approaches for incorporating the PDE in the loss function have been introduced, for instance, in~\cite{weinan2017,yu2018deep,sirignano2018}.
Here, in addition to the PINN loss function, we will consider the loss function of the deep Ritz method~\cite{yu2018deep}, which is based on an equivalent energy minimization problem of the second order elliptic problem. The method is also applicable to other energy minimization problems such as $p$-Laplace problems, contact problems, and elasticity problems.
In particular, the energy minimization problem
\begin{equation*}
	\min_{v \in H^1(\Omega),\;
		v =g \text{ on } \partial \Omega
	}\left( \frac{1}{2} \int_{\Omega} | \nabla v({\bf x})  |^2\, d {\bf x}-\int_{\Omega} f({\bf x}) v({\bf x}) \, d {\bf x}\right),
\end{equation*}
is employed to construct the following practical loss function to train the neural network solution $U({\bf x};\theta)$,
\begin{equation}\label{loss:deepRitz}
	J_{R,M}(\theta)
	:= 
	\frac{w_I}{|X(\Omega)|} \sum_{{\bf x} \in X(\Omega)} \left( \frac{1}{2} |\nabla U({\bf x};\theta)|^2- f({\bf x})U({\bf x};\theta) \right) + \frac{w_B}{|X(\partial \Omega)|} \sum_{{\bf x} \in X(\partial \Omega)} (U({\bf x};\theta)-g({\bf x}))^2.
\end{equation}
Here, again, the boundary condition is enforced with the $L^2$-integral of the error, $U({\bf x};\theta)-g({\bf x})$, and the integrals of the energy term and the boundary condition term are approximated by the Monte--Carlo method; the weight factors $w_I$ and $w_B$ are analogous to those in the PINN loss.

The integral form of the deep Ritz loss function, before approximation via numerical integration, reads
\begin{equation}\label{def:JR}
J_R(\theta):=w_I \int_{\Omega} \left(\frac{1}{2} |\nabla U({\bf x};\theta)|^2 -f({\bf x}) U({\bf x};\theta)\right) \, d{\bf x} + w_B \int_{\partial \Omega} (U({\bf x};\theta)-g({\bf x}))^2 \, ds({\bf x}).
\end{equation}
The subscript $R$ indicates that the loss is formed by the deep Ritz formulation.
In addition, the subscript $M$ in the loss $J_{R,M}(\theta)$ in~\cref{loss:deepRitz} means that the integral in the deep Ritz loss $J_R(\theta)$ in~\cref{def:JR} is approximated by the Monte--Carlo method.

When the boundary condition is implemented as hard constraints using ansatz functions, that is, using the neural network function $\widetilde{U}({\bf x};\theta)$ in~\cref{eq:hard_constraints}, we obtain the integral loss function
$$
	J_{R_I}(\theta)
	:=
	\int_{\Omega} \left( \frac{1}{2} |\nabla \widetilde{U}({\bf x};\theta)|^2-  f({\bf x}) \widetilde{U}({\bf x};\theta)\right) \, d{\bf x}
$$
and train the neural network $\widetilde{U}({\bf x};\theta)$ for the loss $J_{R_I,M}(\theta)$ by approximating the integral in $J_{R_I}$ using the Monte--Carlo method.

It has been numerically studied that for the Poisson model problem, 
the trained solution $U({\bf x};\theta_{P})$ with the PINN loss $J_{P,M}(\theta)$ gives better training results than the trained solution $U({\bf x};\theta_{R})$ with the deep Ritz loss $J_{R,M}(\theta)$; see~\cite{shi2021comparative}. In our work, we will reinvestigate the performance of the two approaches for various hyper parameter settings and report some of our new findings.

As an enhancement to soft enforcement of boundary conditions, an augmented Lagrangian term can be included to the loss function~\cite{Son2023} to obtain,
$$
	L_{P,M}(\theta,\lambda)
	:=
	J_{P,M}(\theta)+\frac{1}{X(\partial \Omega)} \sum_{{\bf x} \in X(\partial \Omega)} (U({\bf x};\theta)-g({\bf x}))\lambda({\bf x})
$$
and
$$
	L_{R,M}(\theta,\lambda)
	:=
	J_{R,M}(\theta)+\frac{1}{X(\partial \Omega)} \sum_{{\bf x} \in X(\partial \Omega)} (U({\bf x};\theta)-g({\bf x}))\lambda({\bf x}),
$$
for the PINN and deep Ritz loss functions, respectively. Here, the boundary condition is enforced as constraints on the neural network solution $U({\bf x};\theta)$
by introducing Lagrange multipliers $\lambda({\bf x})$ for each collocation point ${\bf x}$ in the training sampling set $X(\partial \Omega)$. Hence, $\lambda({\bf x})$ are additional parameters that have to be trained, 
in addition to $\theta$. The use of such an augmented Lagrangian term can improve slow training progress for the boundary loss term and can provide a more accurate trained neural network solution, $U(x;\theta)$. In the augmented Lagrangian approach, the parameters $\theta$ and $\lambda$ are then optimized for the PINN and deep Ritz loss functions in the following sense:
$$
	(\theta_{P},\lambda_{P})
	:=
	\text{arg} \left( \max_{\lambda} \min_{\theta} L_{P,M}(\theta,\lambda)  \right)
	\quad \text{resp.} \quad 
	(\theta_{R},\lambda_{R})
	:=
	\text{arg} \left( \max_{\lambda} \min_{\theta} L_{R,M}(\theta,\lambda)  \right).
$$

We note that the above optimization problems for $\theta$ are non-linear and non-convex while those for $\lambda$ are linear. We thus use the Adam optimization method~\cite{adam} in the gradient update for $\theta$ with a small learning rate $\epsilon$ and a simple gradient update for $\lambda$ with a learning rate $\alpha$, i.e.,
$$
	\lambda=\lambda+ \alpha \nabla_{\lambda} L_{P,M} 
	\quad \text{or} \quad 
	\lambda=\lambda+\alpha \nabla_{\lambda} L_{R,M}.
$$
The learning rate $\alpha$ is often set to a larger value than the learning rate $\epsilon$, as proposed in two-scale update schemes for min-max optimization problems; see~\cite{heusel2017gans,lin2020gradient,chae2023two}. In our numerical experiments, we set $\epsilon=0.001$ for the Adam optimizer and $\alpha=1$ for the gradient ascent update.

We note that the augmented Lagrangian method can be considered as a loss balancing scheme, and in our numerical experiments, we will also conduct comparisons on various loss balancing schemes as listed in~\cref{table:list:para}.

\subsection{Training sampling sets via Gaussian quadrature} \label{sec:gauss}

We recall that, in the loss function of PINN and deep Ritz formulations, the training data sets for $X(\Omega)$ and $X(\partial \Omega)$, and the weight factors $w_I$ and $ w_B$, are the hyper parameters. The training performance and accuracy in the neural network approximation are highly affected by the choice of these hyper parameters.

The Monte--Carlo integration method has a dimension-independent convergence rate and is therefore necessary to beat the curse of dimensionality in high-dimensional domains. In our test problems, the solutions are smooth and the problem domain is a two- or three-dimensional bounded region, and the Monte--Carlo integration method does not take any advantage of such beneficial properties. We note that Gaussian quadrature is recommended for reasonably low-dimensional cases, e.g., in less than five dimensions, and it can also improve the accuracy in the loss computation and the trained solution, see~\cite{mishra2022estimates}. 

Assuming that our model problem is defined in low dimension and has a smooth solution, we propose training sampling sets $X(\Omega)$ and $X(\partial \Omega)$ that are obtained from the Gaussian quadrature; we will employ this in our numerical experiments for the two- and three-dimensional domains. For the two-dimensional case, let the domain $\Omega$ be a rectangle $(a_1,b_1)\times(a_2,b_2)$.
Therefore, we define the following mapping from $(-1,1)$ onto a given interval $(a_k,b_k)$,
$$
	L_k(x)=\frac{a_k(1-x)+b_k(1+x)}{2}.
$$
Then, we can choose $n_G$ Gaussian quadrature points from the interval $(-1,1)$ and corresponding weights $\{ (x_i,w_i) \}_{i=1}^{n_G}$ and transform the quadrature points $\{ x_i \}_{i=1}^{n_G}$ into the points $\{ L_k(x_i) \}_{i=1}^{n_G}$ in the interval $(a_k,b_k)$. In particular, we set
$$
	X_{G}(\Omega)=\{ (L_1(x_i),L_2(x_j))\,:\, \forall i,j=1,\ldots,n_G \}
$$
and similarly
$$
	X_G(\partial \Omega)=\{ (a_1,L_2(x_i)), (b_1,L_2(x_i)), (L_1(x_i),a_2), (L_1(x_i),b_2) \,: \, \forall i=1,\ldots,n_G\}.
$$
Here, for each ${\bf x}=(L_1(x_i),L_2(x_j))$ in $X_G(\Omega)$, we set the associated weight factor $w({\bf x})=\widetilde{w}_i^{(1)} \widetilde{w}_j^{(2)}$. Moreover, for each ${\bf x}=(a_1,L_2(x_i))$ in $X_G(\partial \Omega)$, we set $w({\bf x})=\widetilde{w}_i^{(2)}$, where $\widetilde{w}_\ell^{(k)}$ are defined as the scaled weight factor $\widetilde{w}_\ell^{(k)}:=((b_k-a_k)/2)w_\ell$. The weight factors are chosen analogously for ${\bf x}=(b_1,L_2(x_i)),\;(L_1(x_i),a_2)$, and $(L_1(x_i),b_2)$.

To indicate that the sampling data sets are obtained via Gaussian quadrature, we use the subscript $G$ for the data sets $X_G(\Omega)$ and $X_G(\partial \Omega)$. For the given $n_G$, the number of data points in the set $X_G(\Omega)$ is $n_G^2$ and that in the set $X_G(\partial \Omega)$ is $4n_G$. For the purpose of the comparison, in the Monte--Carlo integration, we also select $n_G^2$ random points from $\Omega$ to form the set $X(\Omega)$ and similarly we form the set $X(\partial \Omega)$ with $4n_G$ randomly chosen points from $\partial \Omega$.

With the training sampling sets $X_G(\Omega)$, $X_G(\partial \Omega)$, we form the loss function
\begin{equation}\label{loss:PNG}
		J_{P,G}(\theta)
		:=w_I \sum_{ {\bf x} \in X_G(\Omega)} \left(\nabla \cdot \nabla U({\bf x};\theta)+f({\bf x}) \right)^2 w({\bf x})\\
		+ w_B \sum_{ {\bf x} \in X_G(\partial \Omega)} \left( U({\bf x};\theta)-g({\bf x}) \right)^2 w({\bf x}),
\end{equation}
in the PINN formulation and 
\begin{equation}\label{loss:DRG}
		J_{R,G}(\theta)
		:=w_I \sum_{ {\bf x} \in X_G(\Omega)} \left( \frac{1}{2}|\nabla U({\bf x};\theta)|^2-f({\bf x})U({\bf x};\theta)\right) w({\bf x})\\
		+ w_B \sum_{ {\bf x} \in X_G(\partial \Omega)} \left( U({\bf x};\theta)-g({\bf x}) \right)^2 w({\bf x}).
\end{equation}
in the deep Ritz formulation. Associated to the above loss functions, we can also form the loss functions with the augmented Lagrangian term,
\begin{equation}\label{loss:LPNG}
	L_{P,G}(\theta,\lambda):=J_{P,G}(\theta)+\sum_{{\bf x} \in X_G(\partial \Omega)} \lambda({\bf x}) (U({\bf x};\theta)-g({\bf x})) w({\bf x})
\end{equation}
and
\begin{equation}\label{loss:LDRG}
	L_{R,G}(\theta,\lambda):=J_{R,G}(\theta)+\sum_{{\bf x} \in X_G(\partial \Omega)} \lambda({\bf x}) (U({\bf x};\theta)-g({\bf x})) w({\bf x}),
\end{equation}
where the boundary condition is enforced as constraints by introducing Lagrange multipliers $\lambda({\bf x})$. In the above, we employ the subscript $G$ to indicate that the integrals in the PINN and deep Ritz loss formulations are approximated by the Gaussian quadrature.

Note that there are also adaptive, residual-based sampling methods, which often improve the performance over simple Monte--Carlo sampling; see, for instance,~\cite{Lu2021DeepXDE:Equations,Nabian2021EfficientSampling,wu2023comprehensive,visser_pacmann_2024}.

\section{Hyperparameters and computation settings}\label{numerics:settings}

In this section, we discuss the hyper parameters under investigation and include details of our computational settings. In~\cref{table:list:para}, we list all hyper parameters considered as well as our recommended hyper parameter choices for the three examples~\eqref{ex1}--\eqref{ex3}; our recommendations will be supported by the numerical results reported in~\cref{numerics:result}. 

\begin{table}
	\begin{center}
		\begin{tabular}{|l|l|c|c|c|}
			\hline
			Hyperparameters & Options & Example 1 & Example 2 & Example 3\\
			\hline
			Loss function & PINN               &   $\checkmark$ & & \\
			& deep Ritz           &  $\checkmark^{*}$  & $\checkmark$ & $\checkmark$ \\ \hline
			Sample sets   & Monte--Carlo method &  $\checkmark$ & & \\
			& Gaussian quadrature   &  $\checkmark^{*}$             & $\checkmark$ & $\checkmark$ \\ \hline
			Loss balance  & weight factor      &  $\checkmark$ & & \\
			& augmented Lagrangian~\cite{Son2023}& $\checkmark$             & $\checkmark$ & $\checkmark$ \\
			& self-adaptive weight~\cite{mcclenny2023self}& $\checkmark$            & & \\
			& inverse Dirichlet~\cite{maddu2022inverse}   & $\checkmark$            & & \\
			& gradient norm~\cite[Algorithm 1 (c)]{wang2023expert}& $\checkmark$           & & 
			 \\ \hline
			\multirow{3}{2cm}{Network architecture enhancements} & boundary condition via ansatz function & $\checkmark$ & $\checkmark$ & $\checkmark$ \\
			& Fourier feature embedding~\cite{tancik_fourier_2020} & $\checkmark$ & & $\checkmark$\\
			& sine activation function & $\checkmark$ & $\checkmark$ & $\checkmark$ \\
			& tanh activation function & $\checkmark$ &  &  \\\hline
			Optimizer & adam~\cite{adam}  & $\checkmark$ &$\checkmark$ & $\checkmark$\\
			& LBFGS~\cite{liu1989limited}  & $\checkmark$ & & \\
			& adam+LBFGS &$\checkmark$ &  & \\
			\hline
		\end{tabular}
	\end{center}
	\caption{List of hyper parameters for numerical study:
		The $\checkmark$ symbol means that the options are recommended for the test examples.
		In Example 1, the superscript in $\checkmark^{*}$ means that
		the deep Ritz formulation and Gaussian quadrature should come in a pair. We also note that a combination, like, PINN formulation and Gaussian quadrature can come in a pair in the above summary of Example 1.}\label{table:list:para}
\end{table}

A summary of the network, sampling set, and optimizer settings that will be used in our computations is listed in~\cref{tab:network:sample:notation}. In particular, as a baseline neural network, we employ a fully connected network with four hidden layers, i.e., $L=4$, and $n$ nodes per each hidden layer with the sine activation function. The number of nodes $n$ per each hidden layer is set differently depending on the complexity of the model problem and the resulting solution. In this context, we compare $\sin$ and $\tanh$ activation functions. Moreover, we test the use of Fourier feature embedding.

As discussed in~\cref{sec:gauss}, we compare Monte--Carlo and Gaussian quadrature schemes to generate training sampling points. For each direction of the problem domain, we choose $n_G$ Gaussian quadrature points to generate the resulting $n_G^2$ interior training sampling points and $4 n_G$ boundary training sampling points. The sum, $n_G^2+4 n_G$, is denoted as $N_t$, the total number of training sampling points. For a fair comparison, in our computations we choose the same number of sampling points $N_t$ in the Monte--Carlo numerical integration.

For the training, we use the Adam optimizer with the learning rate $\epsilon=0.001$ for $\theta$ and the gradient ascent method with the learning rate $\alpha=1$ for the Lagrange multipliers $\lambda$.
We note that, for each training epoch, the parameters $\theta$ and $\lambda$ are updated simultaneously using the Adam optimizer and the gradient ascent method, respectively.
We then train the neural network for a pre-defined maximum number of epochs $T$. 

\begin{figure}[ht!]
	
	\begin{center}
		\includegraphics[width=0.5\textwidth]{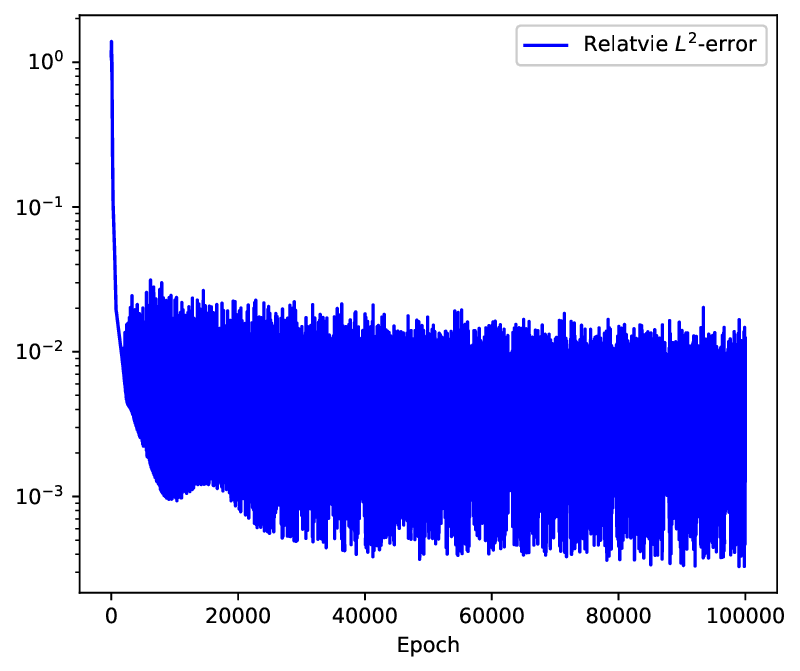}
	\end{center}
	\vskip-.5truecm
	\caption{Relative $L^2$-error history for $U({\bf x};\theta)$ over training epochs for the model solution~\eqref{ex1} with $k=1$: $L_{R,G}$ with $w_I=w_B=1$ is used for the loss function to train the neural network solution $U({\bf x};\theta)$. The error is computed by using a uniform test sample set of $101 \times 101 $ grids over the problem domain.}
	\label{fig:rel_error_hist_smooth}
\end{figure}

As shown in ~\Cref{fig:rel_error_hist_smooth}, the relative $L^2$-error of the neural network during the training can often be smaller than that obtained from the final training epoch. To obtain the trained parameters with a smaller error, we define and use the following error indicators:
\begin{equation*}
	\begin{split}
		E_{P}(\theta)&=\int_{\Omega} (\nabla\cdot\nabla U({\bf x};\theta)+f({\bf x}))^2 \, d{\bf x} + \int_{\partial \Omega} (U({\bf x};\theta)-g({\bf x}))^2 \, ds({\bf x}), \\
		E_{R}(\theta)&=\left\vert \int_{\Omega} \left(|\nabla U({\bf x};\theta)|^2-f({\bf x})U({\bf x};\theta)\right) \, d{\bf x}-\int_{\partial \Omega} \frac{\partial U}{\partial n}({\bf x};\theta) g({\bf x}) \, ds({\bf x}) \right\vert  \\
		& + \int_{\partial \Omega} (U({\bf x};\theta)-g({\bf x}))^2 \, ds({\bf x}),
	\end{split}
\end{equation*}
for the PINN and deep Ritz cases, respectively.
We store the parameters corresponding to the smallest error indicator value observed during the whole training process and use these parameters as the final solution. Unlike the PINN case, the error indicator in the deep Ritz case is set differently from its loss function. We note that the value of the deep Ritz loss function is related to the energy functional and is thus not appropriate for an error indicator. 
The first term in the error indicator $E_R$ is obtained from the  weak formulation of the Poisson problem
by taking the neural network solution $U({\bf x};\theta)$ as a test function. The smaller $E_R$ value
thus indicates that the neural network solution $U({\bf x};\theta)$ is more accurate. In addition, the value $E_R$ can be computed by the first derivatives on $U({\bf x};\theta)$ in contrast to the $E_P$ case where more computation cost is needed for the second derivative calculation.

\begin{table}
	\begin{center}
		\begin{tabular}{|ll|ll|ll|}
			\hline
			\multicolumn{2}{|l|}{Network} & \multicolumn{2}{l|}{Sample points} & \multicolumn{2}{l|}{Optimizer} \\ \hline
			\multicolumn{2}{|l|}{\bf fully connected} & $n_G$: & number of Gaussian & {\bf Adam}: & learning rate $\epsilon=0.001$ \\
			depth: & $L=4$ & & quadrature & Gradient ascent & \multirow{2}{3cm}{$\alpha=1$} \\
			width: & $n$ & $N_t$: & number of total samples & (augm. Lagrange): & \\
			activation: & $\sin$ & & & $T$: & number of training epochs \\
			 &  & & & $E_P$, $E_R$: & error indicators \\
			\hline
		\end{tabular}
	\end{center}
	\caption{Summary of notations for network, sampling points, and optimizer settings.}\label{tab:network:sample:notation}
\end{table}

In our numerical computation, we report the average and the standard deviation of the relative $L^2$-error values for the trained solutions with five different random initializations to show the robustness of our results. We note that, for the augmented Lagrangian approach, we simply initialize all the Lagrange multipliers $\lambda({\bf x})$ by the value $1.0$ and initialize the parameters $\theta$ randomly using a Glorot uniform initializer.
The relative $L^2$-errors are computed by using a test sampling set constructed on a uniform grid of size $101 \times 101$ over the problem domain. 
For the Gaussian quadrature case, for a fixed number of quadrature points $n_G$, the sampling set is also fixed. On the other hand, for the Monte--Carlo numerical integration, the sampling points are randomly initialized with a different random seed in every training run.

Our code has been implemented using the Python JAX library~\cite{bradbury_jax_2018} and the computation is performed on an Intel(R) Xeon(R) Silver 4214R CPU @ 2.40GHz and a Quadro RTX 6000 GPU.

\section{Numerical study on test examples}\label{numerics:result}

In this section, we present the numerical results on comparing the different hyper parameter settings listed in~\cref{table:list:para} for the model problems listed in~\cref{sec:model_problems}. In particular, we first compare different sampling sets in~\cref{sec:sampling}, different loss formulations in~\cref{subsc:loss:study}, network architecture enhancements in~\cref{sec:network_architecture}, loss balancing schemes in~\cref{sec:loss_balancing}, and optimizers in~\cref{sec:optimizers}.

\subsection{Study on sampling sets} \label{sec:sampling}

In this subsection, we test the performance of the PINN and deep Ritz approaches depending on the choice of sampling sets. In our computations, we consider the smooth example in ~\cref{ex1} with $k=1$; see ~\Cref{fig:examples} (left) for the solution. We choose a network with width $n=35$, which leads to a total of $3\,921$ parameters. Moreover, we employ Gaussian quadrature with $n_G=64$, giving $N_t=4\,352$ sampling points, which include $4\,096$ interior and $256$ boundary points to train the network. Therefore, we also randomly select $4\,096$ interior sampling points and $256$ boundary sampling points in the case of Monte--Carlo integration. We employ the loss functions with or without the augmented Lagrangian term  and  train network parameters $\theta$ and the Lagrange multipliers $\lambda$ for $T=100\,000$ epochs. For the sake of clarity, we summarize the loss function notations depending on the loss function formulations and the integration schemes in~\cref{tab:loss:sample}.

\begin{table}[t]
	\begin{center}
		\begin{tabular}{|c|c|c|}
			\hline
			& Monte--Carlo integration & Gaussian quadrature\\ \hline
			PINN        & $J_{P,M}$   & $J_{P,G}$ \\
			deep Ritz    & $J_{R,M}$   & $J_{R,G}$ \\
			PINN-AL     & $L_{P,M}$   & $L_{P,G}$ \\
			deep Ritz-AL & $L_{R,M}$   & $L_{R,G}$ \\
			\hline
		\end{tabular}
	\end{center}
	\caption{Notations for various loss formulations and sampling sets: PINN (standard PINN loss), deep Ritz (standard deep Ritz loss), PINN-AL (PINN loss with the augmented Lagrangian term), deep Ritz-AL (deep Ritz loss with the augmented Lagrangian term), Monte--Carlo (Monte--Carlo numerical integration), and Gaussian quadrature; cf.~the discussion in~\cref{sec:nn_approximation,sec:gauss}.}
	\label{tab:loss:sample}
\end{table}

\begin{figure}[ht!]
	
	\begin{center}
		\includegraphics[width=1.0\textwidth]{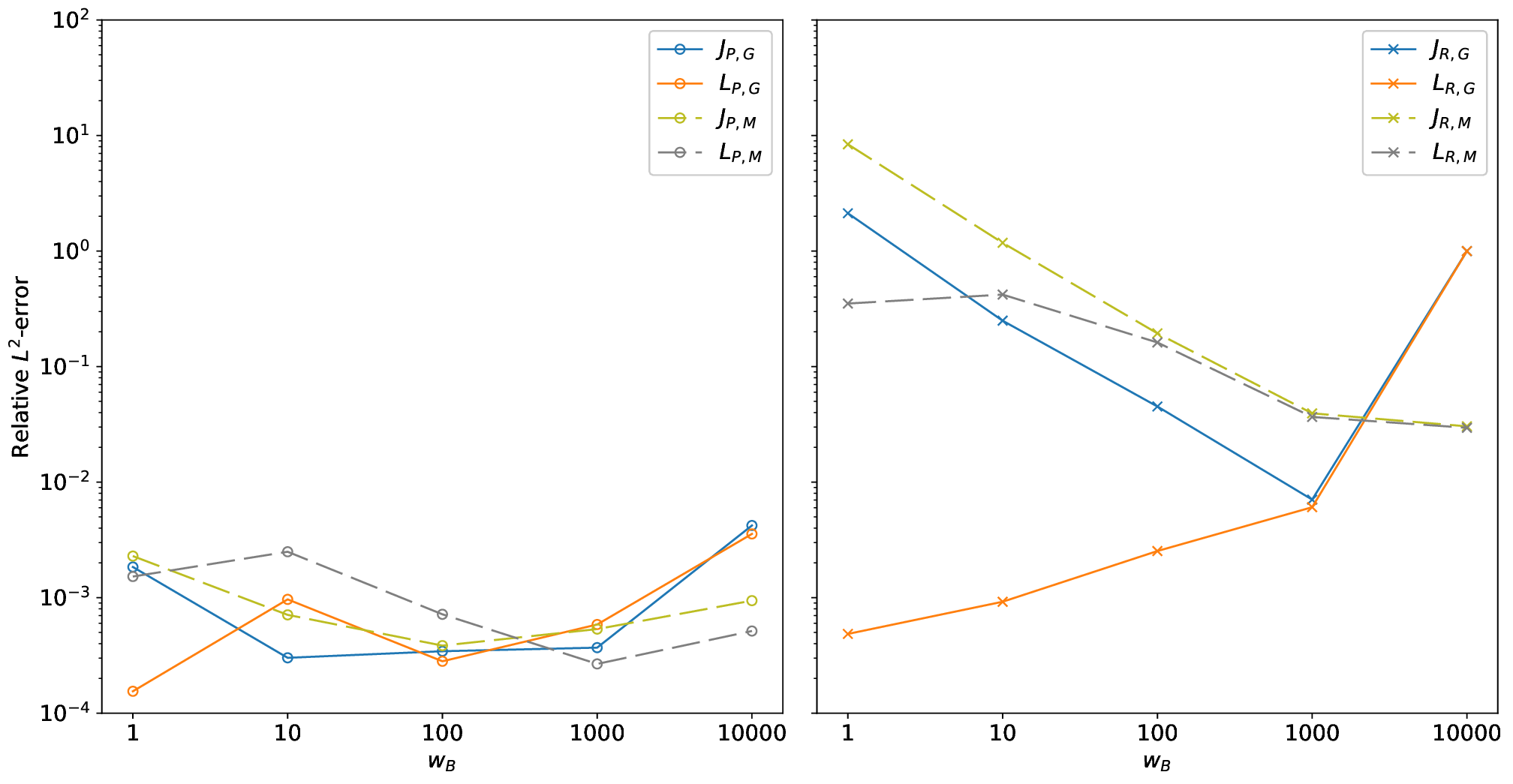}
	\end{center}
	\vskip-.5truecm
	\caption{Study on sampling sets for the example in \eqref{ex1} with $k=1$: the average of the relative $L^2$-errors over various $w_B$ choices for the four loss formulations (left: PINN, right: deep Ritz) depending on the sampling approach (solid line: Gaussian quadrature, dashed line: Monte Carlo).}
	
	\label{fig:plot-Table4}
\end{figure}

In~\cref{tab:sample}, the relative $L^2$-errors of the neural network approximation to the exact solution are reported. For the weight factors in the loss function, we set $w_I=1$ and various values for the weight factor $w_B = 1, 10, 100, 1\,000, 10\,000$. For the standard PINN loss $J_P$, there is no significant difference in the obtained results depending on the integration methods, while the deep Ritz loss case $J_R$ can benefit from using Gaussian quadrature instead of Monte--Carlo integration. In addition, the standard PINN loss yields similar errors with and without the augmented Lagrangian approach depending on the integration methods.
However, for the deep Ritz formulation with Gaussian quadrature, the use of the augmented Lagrangian term significantly improves the results, i.e., $L_{R,G}$ performs much better than $J_{R,G}$.
{In Figure~\ref{fig:plot-Table4}, {we plot} the error {while varying} the weight factor $w_B${. The results} clearly show the benefit of using the Gaussian quadrature sample set for both the PINN and deep Ritz { formulation} with the augmented {Lagrangian} term. }

\begin{table}[ht!]
		\caption{Study on sampling sets for the example in \eqref{ex1} with $k=1$: the average of the relative $L^2$-errors and their standard deviation (inside the parenthesis).
		} \label{tab:sample}
		{\normalsize \renewcommand{\arraystretch}{1.0}
			\begin{center}
				\vskip-.3truecm
				\begin{tabular}{ccccccc}
					
					\hline\hline
					& $w_B$ & 1  & 10 & 100 & 1000 & 10000 \\
					\Xhline{3\arrayrulewidth}
					\multirow{8}{2cm}{Gaussian quadrature} & $J_{P,G}$    		   & 1.843e-03      & {\bf3.016e-04} & 3.439e-04 & 3.694e-04      & 4.213e-03 \\
					&						   & (1.36e-03)     &   (4.80e-05)   & (8.00e-05)& (1.54e-04)     & (1.62e-03)\\
					& $J_{R,G}$    		   & 2.124e-00      &    2.496e-01   & 4.516e-02 & {\bf7.060e-03} & 9.994e-01 \\
					&						   & (4.06e-03)     &   (3.71e-03)   & (4.13e-03)& (9.74e-04)     & (1.31e-03)\\
					& $L_{P,G}$    		   & {\bf1.548e-04}      & 9.666e-04 & 2.814e-04 & 5.858e-04      & 3.563e-03 \\
					&						   & (4.94e-05)     &   (1.91e-05)   & (5.47e-05)& (1.06e-04)     & (8.71e-04)\\
					& $L_{R,G}$    		   & {\bf4.854e-04} &    9.220e-04   & 2.533e-03 & 6.065e-03      & 1.000e-00 \\
					&						   & (2.35e-04)     &   (3.62e-04)   & (6.11e-04)& (9.53e-04)     & (2.23e-04)\\
					\Xhline{0.8\arrayrulewidth}
					\multirow{8}{2cm}{Monte--Carlo integration}	  & $J_{P,M}$    	       & 2.291e-03 & 7.100e-04 & {\bf3.846e-04} &   5.347e-04    &    9.424e-04   \\
					&						   & (5.62e-04)& (1.65e-04)&   (1.50e-04)   &   (4.88e-04)   &   (7.48e-04)   \\
					& $J_{R,M}$    	 	   & 8.401e-00 & 1.178e-00 &   1.940e-01    &   3.941e-02    & {\bf3.040e-02} \\
					&						   & (1.81e-01)& (1.16e-01)&   (1.02e-01)   &   (1.12e-02)   &   (6.72e-03)   \\
					& $L_{P,M}$     		   & 1.524e-03 & 2.497e-03 &   7.179e-04    & {\bf2.670e-04} &    5.145e-04   \\
					&						   & (3.97e-04)& (1.80e-04)&   (1.80e-04)   &   (8.58e-05)   &   (1.58e-04)   \\
					& $L_{R,M}$     		   & 3.514e-01 & 4.200e-01 &   1.617e-01    &   3.666e-02    & {\bf2.953e-02} \\
					&						   & (4.17e-02)& (1.64e-01)&   (3.34e-02)   &   (1.12e-02)   &   (7.11e-03)   \\
					\Xhline{3\arrayrulewidth}
				\end{tabular}
			\end{center}
		}
		\vskip-.2truecm
\end{table}

\begin{figure}[ht!]
	\begin{center}
		\includegraphics[width=0.45\textwidth]{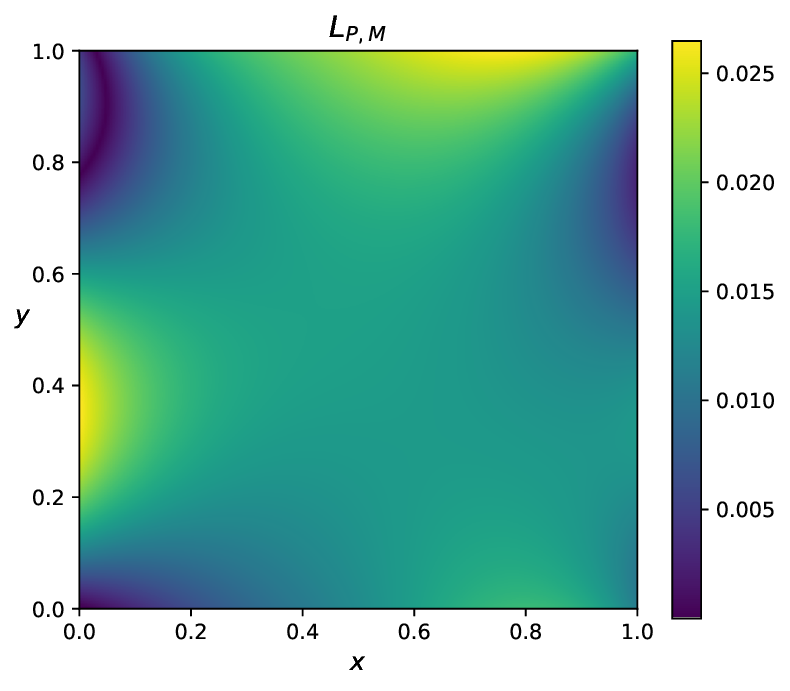} \qquad
		\includegraphics[width=0.45\textwidth]{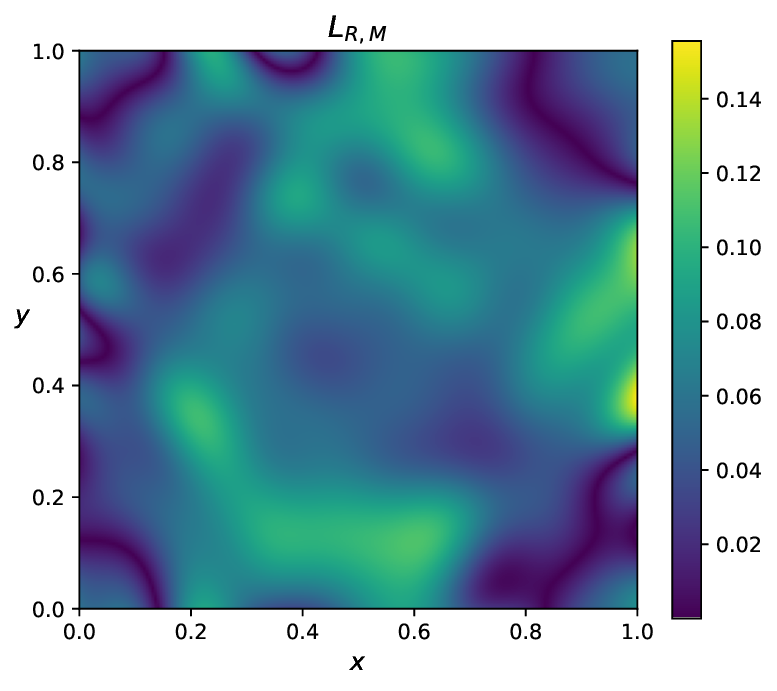} \\
		\includegraphics[width=0.45\textwidth]{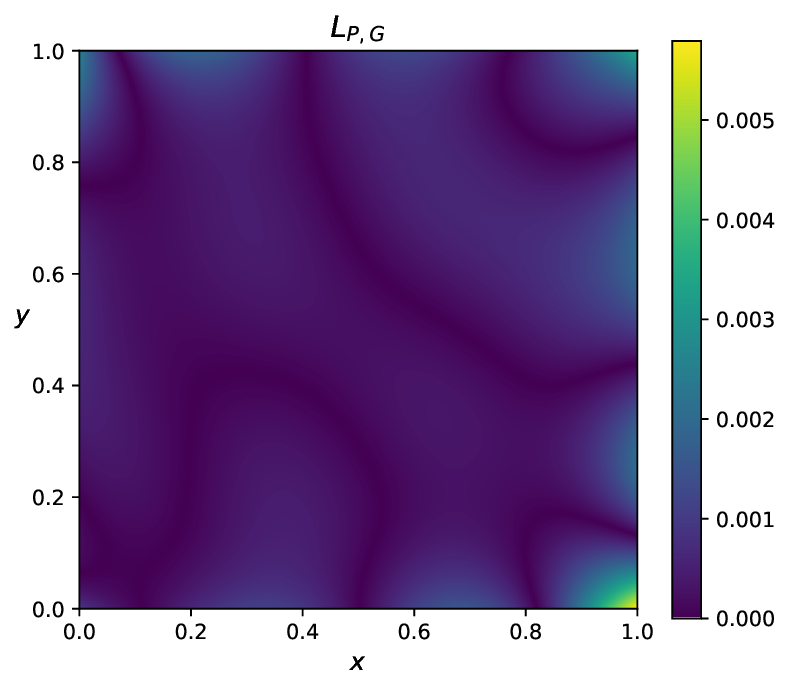} \qquad
		\includegraphics[width=0.45\textwidth]{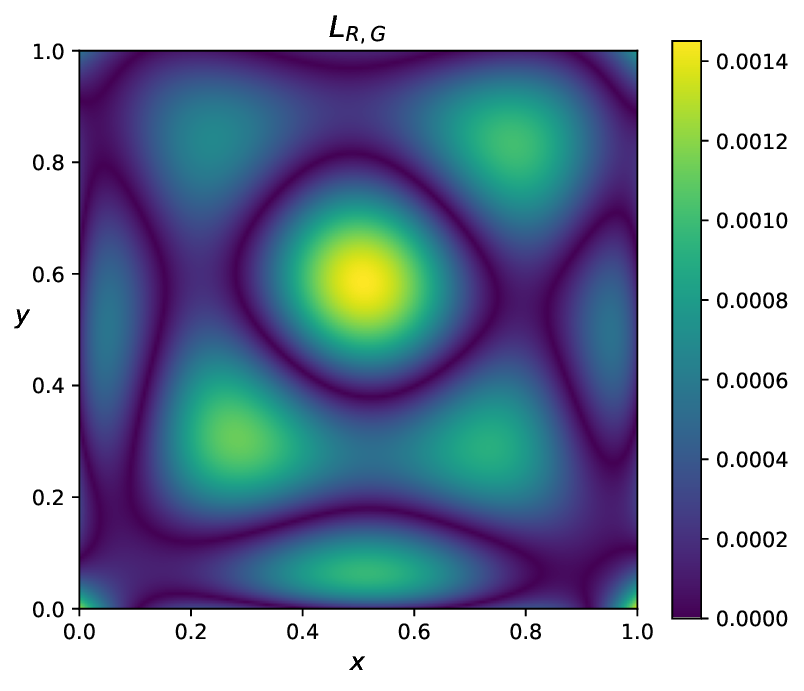}
	\end{center}
	\vskip-.7truecm
	\caption{Smooth example in~\eqref{ex1} with $k=1$: absolute error plots for the PINN loss $L_{P,M}$ (top, left) and the deep Ritz loss $L_{R,M}$ (top, right) using Monte--Carlo integration as well as the PINN loss $L_{P,G}$ (bottom, left) and the deep Ritz loss $L_{R,G}$ (bottom, right) using Gaussian quadrature.}
	\label{fig:smooth_error}
\end{figure}

\Cref{fig:smooth_error} shows plots of the absolute error against the exact solution of the neural network solutions that have been trained using the PINN and deep Ritz approaches with $w_I=w_B=1$ and the augmented Lagrangian term, using either Monte--Carlo integration or Gaussian quadrature. The use of Gaussian quadrature yields smaller errors than Monte--Carlo integration for both PINN and deep Ritz cases. When using the Gaussian quadrature, both $L_{P,G}$ and $L_{R,G}$ loss cases show similar relative $L^2$-error values, but the absolute error plot for the PINN loss case, i.e., $L_{P,G}$, shows higher error values near the boundary, while that for the deep Ritz loss case, i.e., $L_{R,G}$, shows evenly distributed error values over the problem domain.

For the two other test examples, we observed a similar behavior depending on the training sampling sets. In the following experiments, we thus restrict ourselves to using Gaussian quadrature, if not mentioned otherwise. We note again that the deep Ritz approach with augmented Lagrangian term gives much more accurate results when Gaussian quadrature sampling points are used. In the next section, we will see that this combination outperforms PINNs for some more challenging model problems, like multi-frequency solutions in~\cref{ex2} or high-contrast and oscillatory solutions in~\cref{ex3}.

\subsection{Study on loss formulations}\label{subsc:loss:study}

In this section, we compare various loss formulations for the three model problems. In this context, we observe that the deep Ritz formulation with the augmented Lagrangian loss outperforms other loss formulations for the multi-component oscillatory model problem and high-contrast oscillatory model problem in~\cref{ex2} and~\cref{ex3}, respectively. In this context, we choose $w_I=1$ and various values for $w_B =10^k$, for $k=0,1,\ldots,4$, for loss functions without augmented Lagrangian. Since the augmented Lagrangian term is included to deal with the loss balance for the boundary condition, we simply set the weight factors $w_I=1$ and $w_B=1$, in these cases. 

For the test example~\eqref{ex1}, we again consider the model solution with $k=1$ and employ a fully connected neural network with width $n=35$, a training sampling set with $n_G=64$ Gaussian quadrature points, and $T=100\,000$ training epochs. For the test example in~\eqref{ex2}, we consider the model solution with $N=6$ and a fully connected neural network with width $n=100$, 
a training sampling set with $n_G=256$ Gaussian quadrature points, and $T=1\,000\,000$ training epochs; a relatively large neural network is required due to the high complexity of the solution for large values of $N$. 
Finally, for the test example in~\cref{ex3}, we consider the model solution with $A=100$ and $\varepsilon=0.01$ and the same network architecture, sampling set, and number of training epochs as in first test problem given in~\eqref{ex1}. The average of the relative $L^2$-error values and their standard deviation are listed in~\cref{tab:loss}. In addition, for comparison purposes the average computation time is reported for the test example in~\eqref{ex2} with the loss formulations $L_{P,G}$ and $L_{R,G}$, that is, for the PINN and deep Ritz methods with augmented Lagrangian loss and Gaussian quadrature.

We observe that, for Example~\eqref{ex1}, the standard PINN with a larger weight factor $w_B$ performs well, while it shows much larger errors and computation time for a more challenging Example~\eqref{ex2}. 
In contrast to PINNS, the use of a larger weight factor in the deep Ritz formulation, i.e., $J_{R,G}$ with larger weight factors $w_B$, does not help to improve the accuracy. Here, the inclusion of the augmented Lagrangian term is more important for the accuracy of the trained model. As mentioned before, this combination, performs best for more challenging examples, like examples in~\cref{ex2,ex3}, that are difficult to train using the considered PINN formulation; see also the discussion in~\cite{NTK-theory,dolean2024multilevel}.
We also would like to stress the reduced computing time for the deep Ritz compared with the PINN formulation. This is since the loss evaluation only requires the computation of the first derivative on the neural network; cf.~the computation time results in~\cref{tab:loss} for $L_{P,G}$ and $L_{R,G}$ of Example~\eqref{ex2}.

\begin{table}[ht!]
	
	\caption{Loss formulation study on test examples: average and standard deviation of the relative $L^2$-errors (in parenthesis), {the {red colored numbers} indicate the average computation time of $L_{P,G}$ and $L_{R,G}$ with $w_B=1$.} }\label{tab:loss}
	{\normalsize \renewcommand{\arraystretch}{1.0}
		\begin{center}
			\vskip-.3truecm
			\begin{tabular}{cccccccc}
				
				\hline\hline
				& $w_B$ & \multicolumn{2}{c}{1}  & 10 & 100 & 1000 & 10000 \\
				\Xhline{3\arrayrulewidth}
				\multirow{8}{*}{Example~\eqref{ex1}}
				& $J_{P,G}$ & 1.843e-03 &   & {\bf3.016e-04} &    3.439e-04  &   3.694e-04    & 4.213e-03  \\
				&           & (1.36e-03) &  &    (4.80e-05)  &    (8.00e-05) &   (1.54e-04)   & (1.62e-03) \\
				& $J_{R,G}$ & 2.124e-00  &  &    2.496e-01   &    4.516e-02  & {\bf7.060e-03} & 9.994e-01  \\
				&           & (4.06e-03) &  &    (3.71e-03)  &    (4.13e-03) &   (9.74e-04)   & (1.31e-03) \\
				& $L_{P,G}$ & \bf 1.548e-04 & \bf \color{red} 440s & & & & \\
				&           & (4.94e-05) &  &   &  &  & \\
				& $L_{R,G}$ & \bf 4.854e-04 & \bf \color{red} 360s & & & & \\
				&           & (2.35e-04) &  &   &  &  & \\
				\Xhline{0.8\arrayrulewidth}
				\multirow{8}{*}{Example~\eqref{ex2}}
				& $J_{P,G}$ & 9.236e-01  &  & 6.014e-01    & 3.614e-01    & {\bf1.267e-01}   & 1.660e-01 \\
				&           & (3.54e-02) &  & (1.90e-01) & (9.20e-02) & (5.56e-02) & (7.39e-02) \\
				& $J_{R,G}$ & 2.103e-00  &  & 4.533e-01    & {\bf1.753e-01}   & 3.093e-01   & 1.000e-00\\
				&           & (2.61e-02) &  & (5.55e-02) & (1.30e-01) & (3.73e-01) & (1.55e-05) \\
				& $L_{P,G}$ & \bf 4.078e-01 & \bf \color{red} 41\,000s & & & & \\
				&           & (1.19e-01) &  &   &  &  & \\
				& $L_{R,G}$ & \bf 6.540e-03 & \bf \color{red} 12\,000s & & & & \\
				&           & (1.44e-03) &  &   &  &  & \\
				\Xhline{0.8\arrayrulewidth}
				\multirow{8}{*}{Example~\eqref{ex3}}
				& $J_{P,G}$ & 2.306e-01  &    & 9.111e-02    & 2.303e-02   & {\bf1.130e-02}   & 1.000e-00 \\
				&           & (3.36e-02) &    & (3.09e-02)   & (7.92e-03)  & (8.97e-03)       & (1.78e-05) \\
				& $J_{R,G}$ & {\bf2.575e-02}  &  & 4.128e-02    & 4.152e-01   & 1.000e-00   & 1.000e-00 \\
				&           & (1.51e-03)  &   & (1.30e-02)  & (4.78e-01) & (0.00e-00) & (0.00e-00) \\
				& $L_{P,G}$ & \bf 1.845e-02 & \bf \color{red} 440s & & & & \\
				&           & (2.00e-03) &    &   &  &  & \\
				& $L_{R,G}$ & \bf 1.407e-03 & \bf \color{red} 360s & & & & \\
				&           & (2.56e-04)  &   &   &  &  & \\
				\Xhline{3\arrayrulewidth}
			\end{tabular}
		\end{center}
	}
	\vskip-.2truecm
	
\end{table}

\begin{figure}[ht!]
	\begin{center}
		\includegraphics[width=0.32\textwidth]{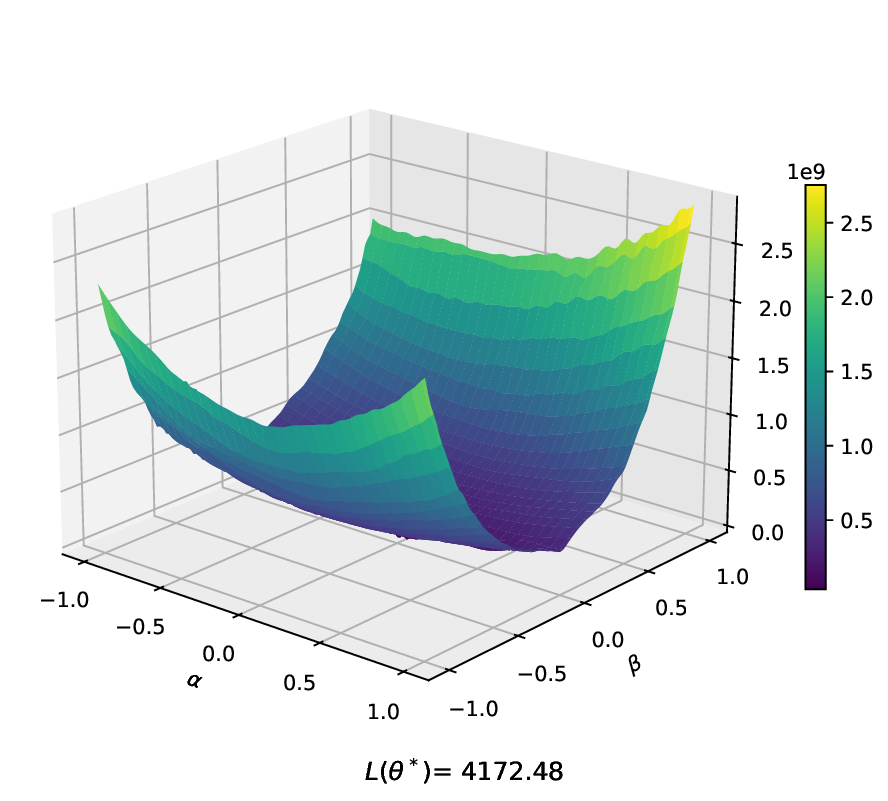}\hfill
		\includegraphics[width=0.32\textwidth]{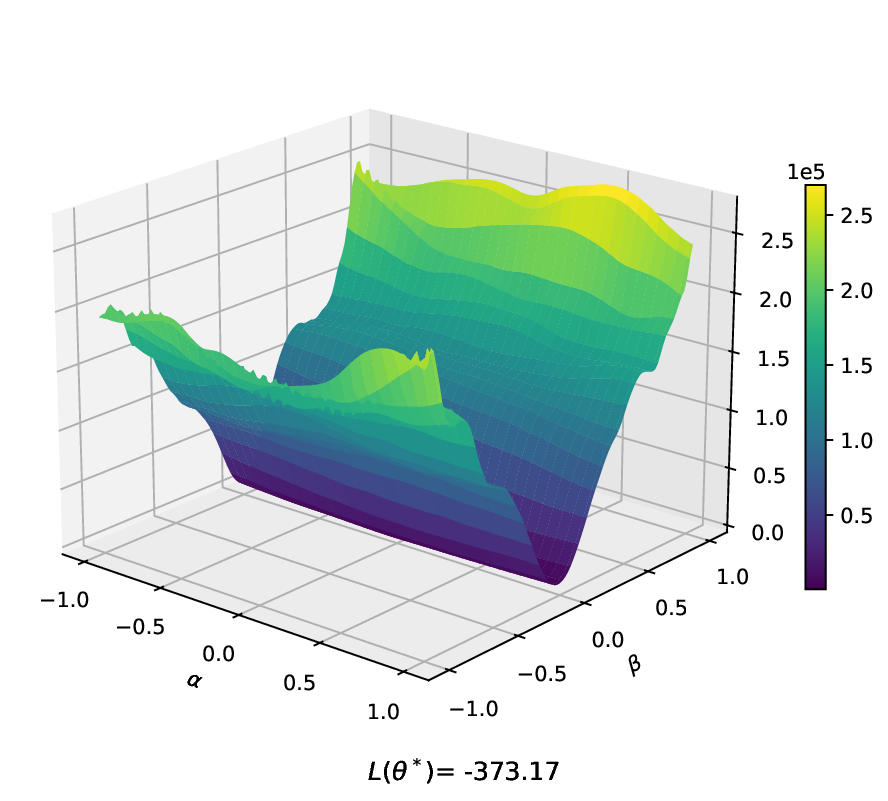}\hfill
		\includegraphics[width=0.32\textwidth]{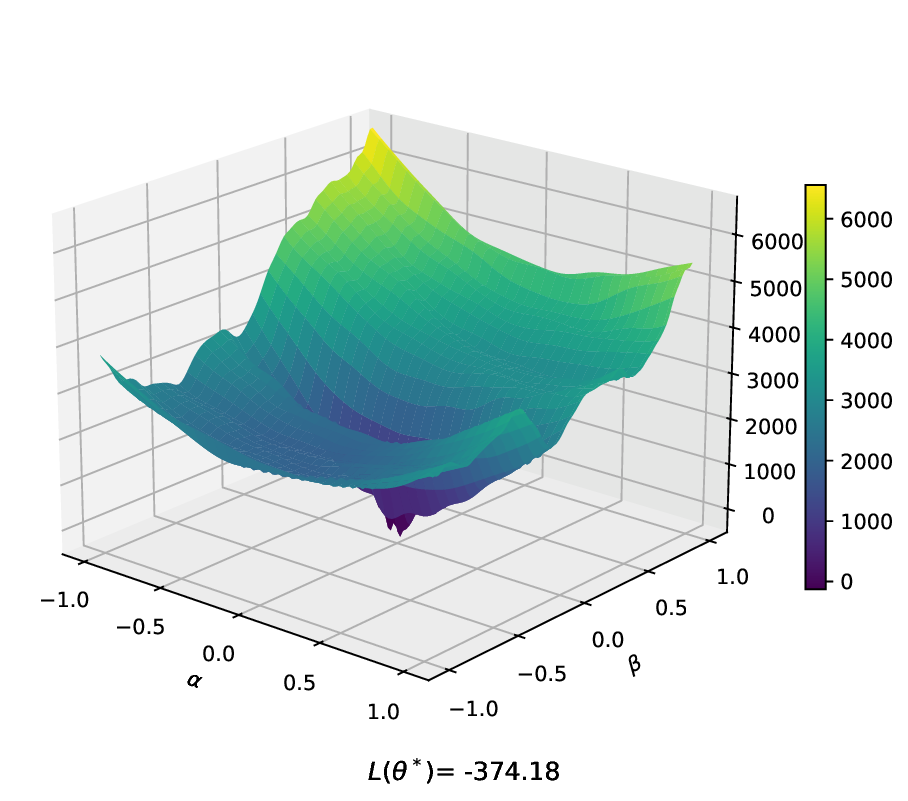}\\
		\includegraphics[width=0.32\textwidth]{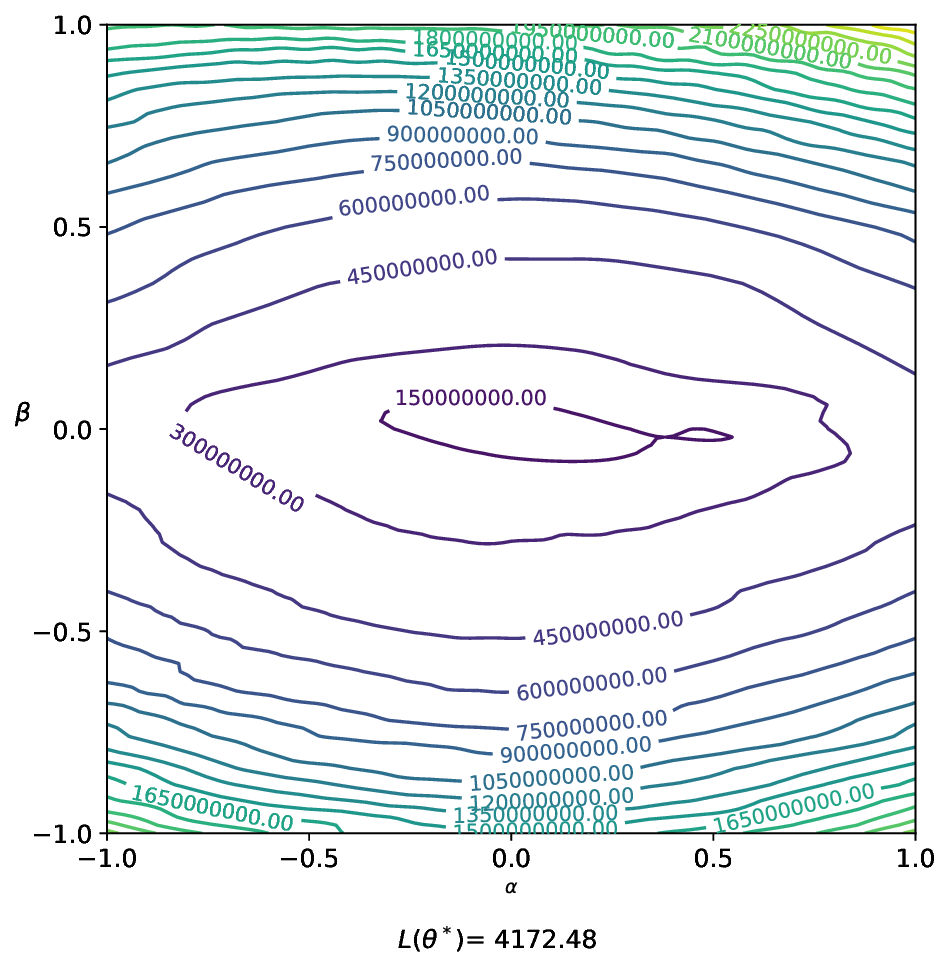}\hfill
		\includegraphics[width=0.32\textwidth]{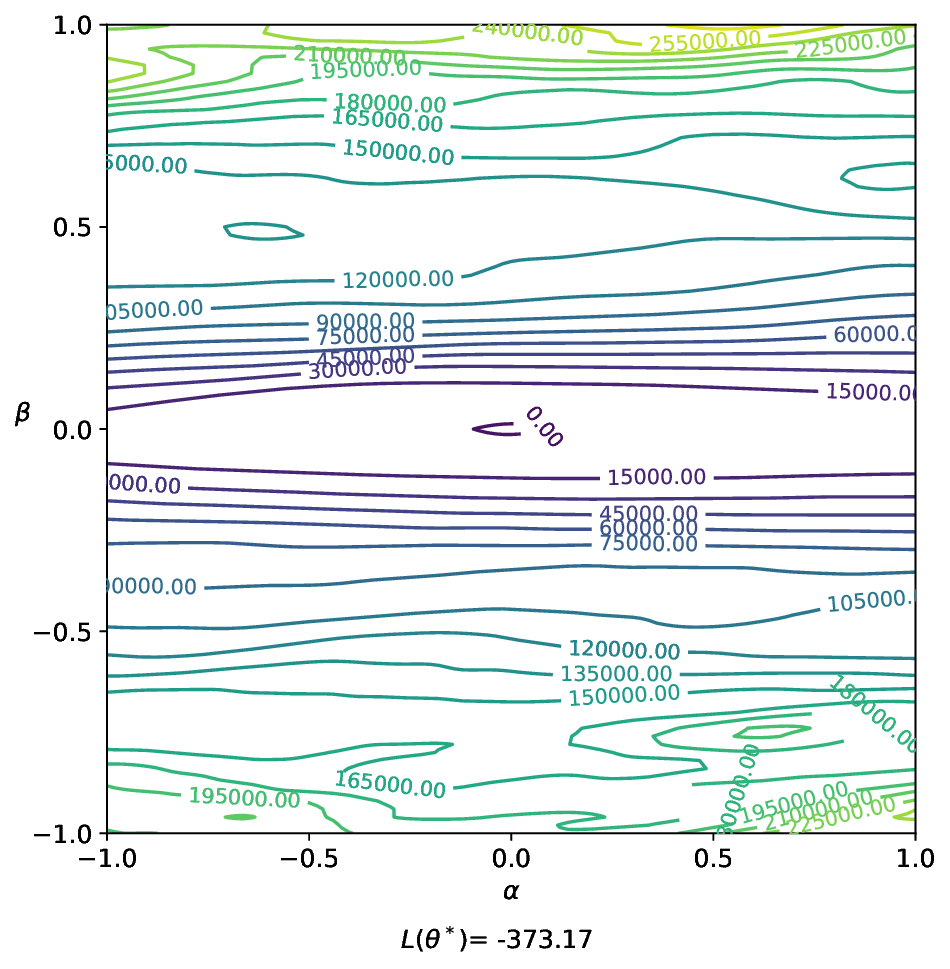}\hfill
		\includegraphics[width=0.32\textwidth]{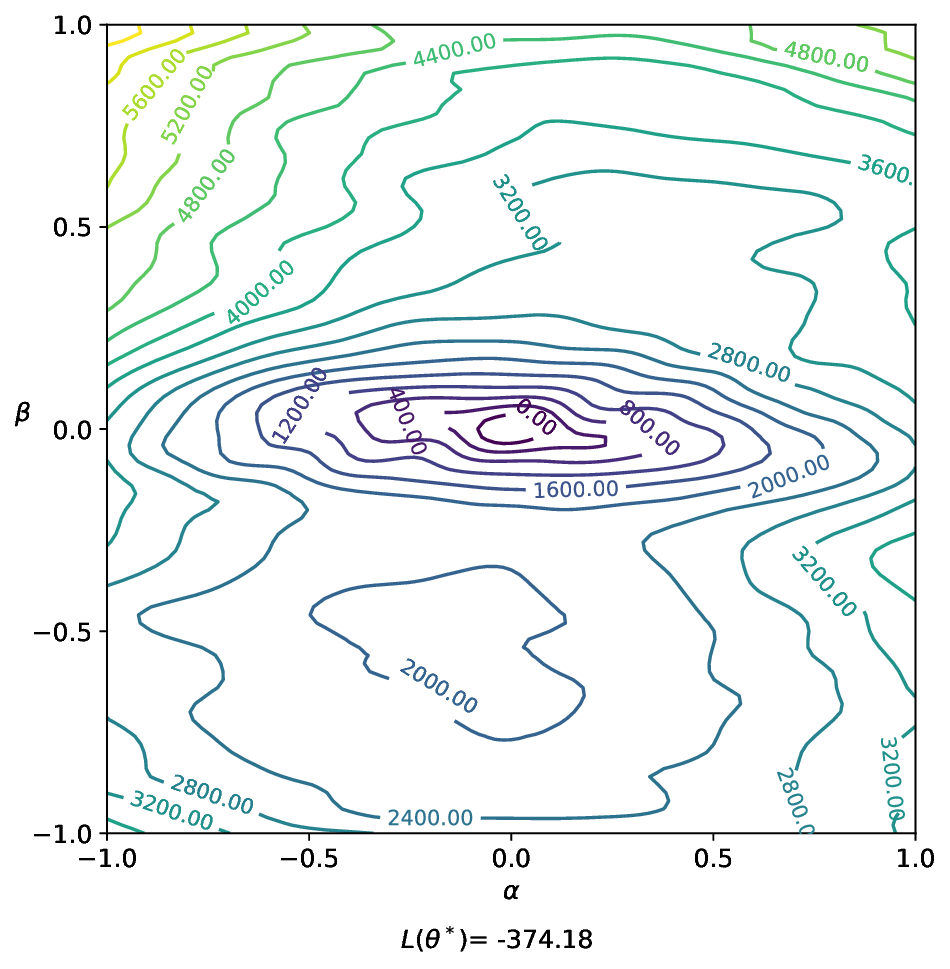}\\
        \includegraphics[width=0.32\textwidth]{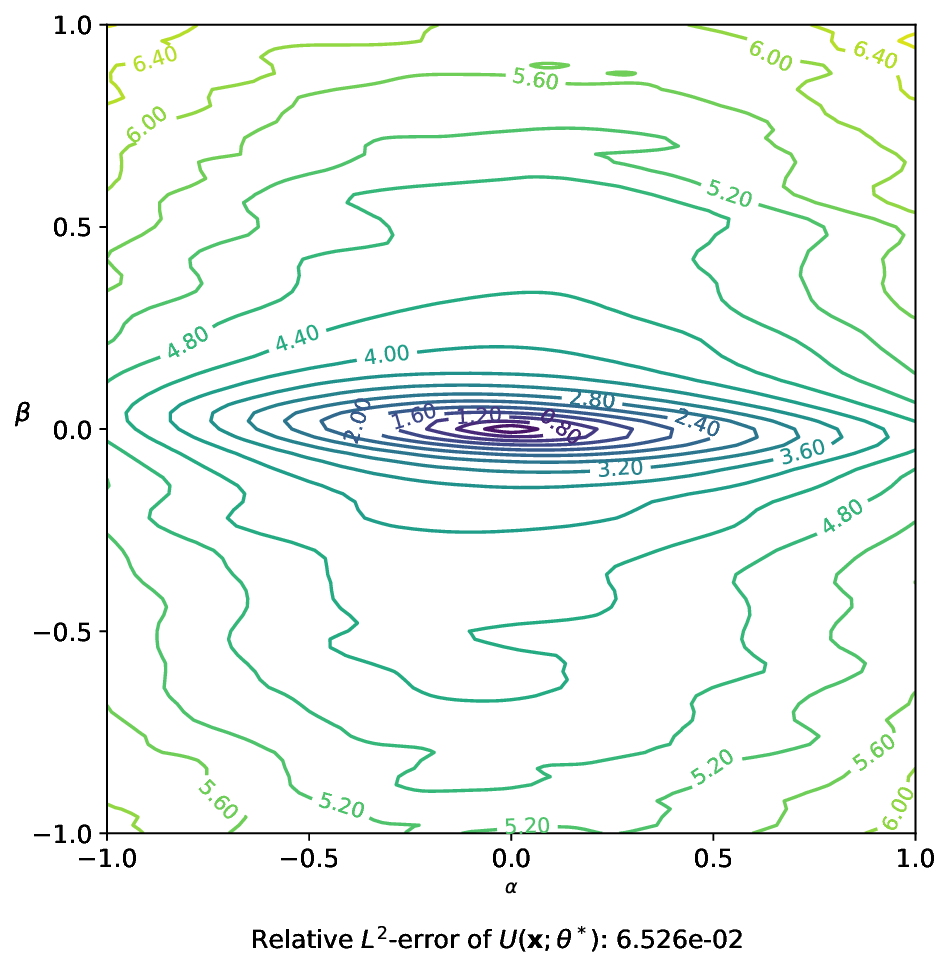}\hfill
		\includegraphics[width=0.32\textwidth]{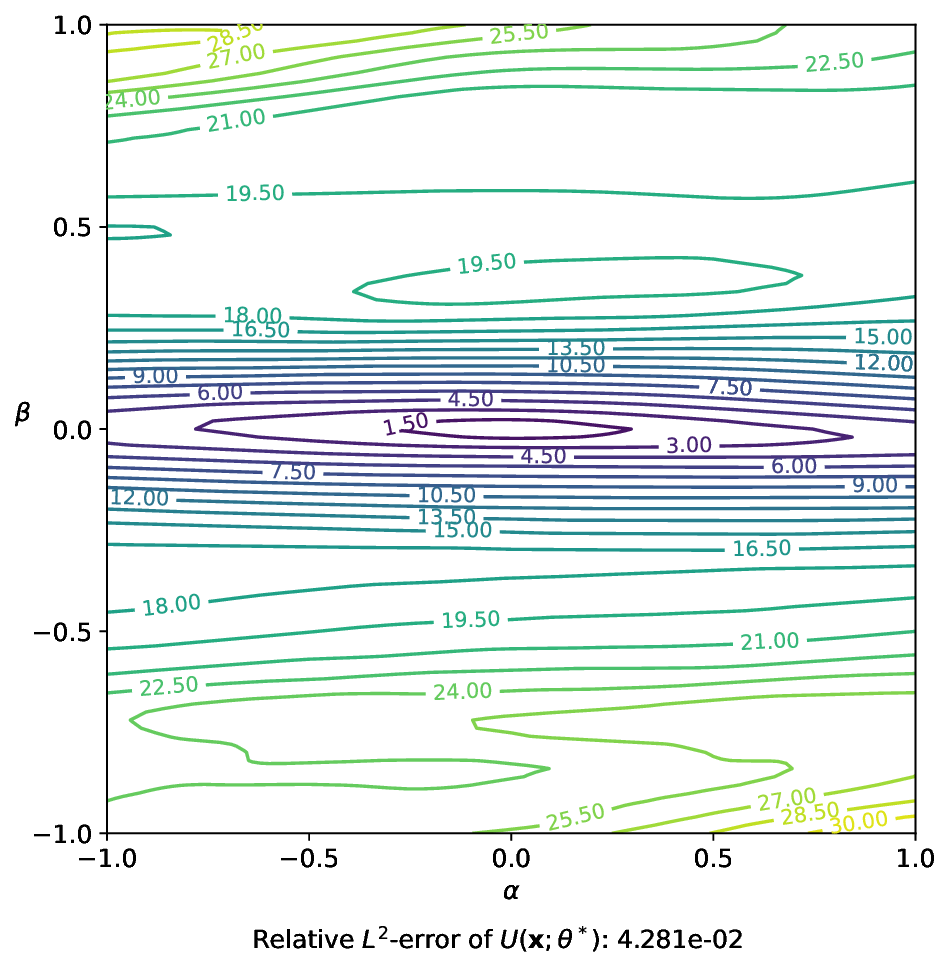}\hfill
		\includegraphics[width=0.32\textwidth]{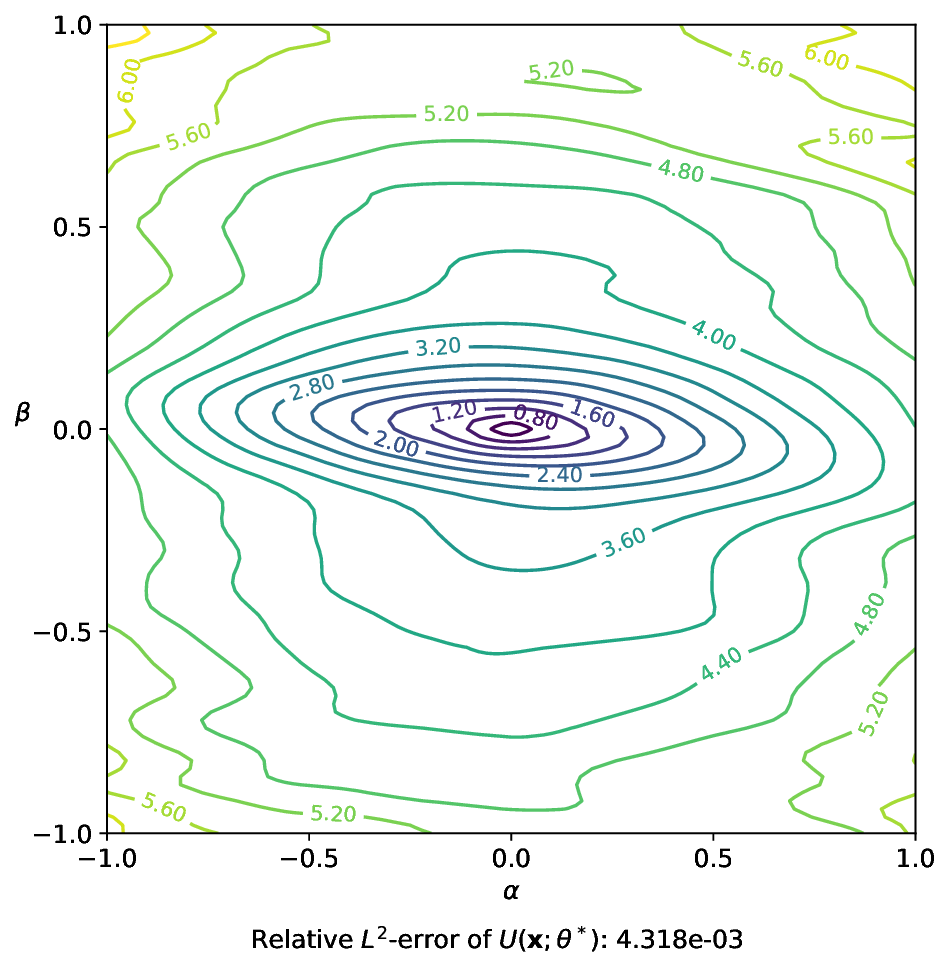}
	\end{center}
	\vskip-.7truecm
	\caption{Loss landscape of Example~\ref{ex2}-Surface plot (top) and
    Contour plot (middle), and the relative $L^2$-error for the corresponding $U({\bf x};\widetilde{\theta})$ in Contour plot (bottom): $J_{P,G}$ (left), $J_{R,G}$ (middle), and $L_{R,G}$ (right)}
	\label{fig:loss-landscape}
\end{figure}

To analyze the effectiveness of the loss formulation $L_{R,G}$ for Example~\eqref{ex2},
we compare the loss landscape of the three loss formulations $J_{P,G}$, $J_{R,G}$, and $L_{R,G}$ by using the visualization method proposed in \cite{li2018visualizing}. For the corresponding trained parameter $\theta^*$ to each loss formulation, we compute the loss function value for the perturbed parameter $\widetilde{\theta}$ obtained from the layer and direction normalization, i.e.,
$$L(\widetilde{\theta}):=L(\theta^*+\alpha 
{\bf \zeta} + \beta {\bf \gamma}).$$
{In the above, the direction vector ${\bf \zeta}$ consists of $({\bf \zeta}_i)_i$
and the $i$-the layer vector ${\bf \zeta}_i$ is computed by}
$${\bf \zeta}_i=\frac{{\bf d}_i}{ \| {\bf d}_i \| }
\| \theta^*_i \|,$$
where $\theta_i^*$ denotes the $i$-th layer parameter of $\theta^*$, ${\bf d}_i$ denotes a random Gaussian direction vector to the parameter $\theta_i^*$, 
and $\| \cdot \|$ denotes the Frobenius norm. The direction ${\bf \gamma}$ is obtained 
analogously with a different random direction.
In \Cref{fig:loss-landscape}, we present the plots of the loss landscape to the three different loss formulations.
We can see that $L_{R,G}$ loss formulation gives a better landscape near the trained $\theta^*$, in the sense that the loss landscape is more uniform in the two different directions and less stretched in one specific direction. This may lead to better convergence of the Adam optimizer.

In this context, recall that the directions are chosen randomly, so this observation may not hold for all directions. However, the observation is in alignment with our observation that the solution is better for the $L_{R,G}$ loss formulation.
In addition, the contour plot of the relative $L^2$ error value for $U({\bf x};\widetilde{\theta})$ shows a similar behavior to the loss landscape plot, {which also indicates that the $L_{R,G}$ loss { yields} better trained solutions compared to the other two loss formulations.}

\subsection{Network architecture enhancements} \label{sec:network_architecture}

In this subsection, we study the training performance and accuracy depending on certain neural network architecture enhancements, i.e., hard enforcement of boundary conditions via an ansatz function, the choice of the activation function, and the inclusion of Fourier feature embedding~\cite{tancik_fourier_2020}.
We will consider the three test examples in~\cref{ex1,ex2,ex3} with $k=4$, $N=4$, and $A=100$ and $\varepsilon=0.01$, respectively. For all the examples, we employ a network with width $n=35$, the same sampling sets with $n_G=64$,
and a total number of $T=100\,000$ training epochs.

\paragraph{Ansatz function} We first study the hard enforcement of boundary conditions via an ansatz function. In our computation, we only consider the zero boundary condition and we simply set $A({\bf x} )=0$ and study various choices for $G( \bf{x} )$ in \cref{eq:hard_constraints}. For the ansatz function $G({\bf x})$, we test and compare the following options:
\begin{equation}
	G({\bf x} )\in\{ x(1-x)y(1-y), \; \sin(\pi x)\sin(\pi y),\;  \sin(4 \pi x)\sin(4 \pi y),\; \sin(8 \pi x)\sin(8 \pi y)\}.
\end{equation}

In~\cref{tab:ansatz}, we report the average and standard deviation of the relative $L^2$-error values for the different choices of the ansatz function $G({\bf x})$ as well as, for the sake of comparison, those obtained from the deep Ritz formulation with the augmented Lagrangian for comparison. We note that the choice $G({\bf x} )=\sin (4 \pi x)\sin( 4 \pi y)$ is identical to the exact solution in Example~\eqref{ex1} and we thus obtained very accurate trained solution with this particular choice. For the examples in~\cref{ex2,ex3}, with the choice $G( {\bf x} )=\sin(\pi x)\sin(\pi y)$ we obtained the best results, with a similar accuracy as those obtained with the $L_{R,G}$ formulation. However, we observe that the ansatz function has to be chosen with care. For instance, with the choice $G({\bf x} )=\sin(4 \pi x)\sin( 4 \pi y)$ the accuracy deteriorates for Examples~\eqref{ex2} and~\eqref{ex3}.

\begin{table}[ht!]
	\caption{Study of hard enforcement of boundary conditions on the test examples: average and standard deviation of the relative $L^2$-errors (in parenthesis); the results obtained from $L_{R,G}$ (boldface and colored in red) are listed for the comparison.} \label{tab:ansatz}
	{\normalsize \renewcommand{\arraystretch}{1.0}
		\begin{center}
			\vskip-.3truecm
			\begin{tabular}{cccccc}
				
				\hline\hline
				& $G({\bf x})$ & $x(1-x)y(1-y)$ & $\sin(\pi x)\sin(\pi y)$ & $\sin(4 \pi x)\sin(4 \pi y)$
				& $\sin(8 \pi x)\sin(8 \pi y)$\\
				\Xhline{3\arrayrulewidth}
				\multirow{6}{*}{Example~\eqref{ex1}}
				& ${J}_{P_I,G}$ & 2.432e-04    & 2.195e-04   & {\bf1.300e-06} & 9.962e-01  \\
				&               & (4.14e-04)   & (7.56e-05) &(5.16e-07)&(4.71e-03) \\
				& ${J}_{R_I,G}$ & 1.314e-03    & 2.885e-03  & {\bf9.460e-05}& 9.606e-01  \\
				&               & (7.33e-03)   & (4.80e-03) &(4.26e-05)&(5.48e-03) \\
				& $L_{R,G}$     & {\bf \color{red} 6.400e-04}     &   &  & \\
				&               & {\bf \color{red} (1.36e-04)}   &   &  & \\
				\Xhline{0.8\arrayrulewidth}
				\multirow{6}{*}{Example~\eqref{ex2}}
				& ${J}_{P_I,G}$ & 2.417e-01    & {\bf1.177e-02}  & 8.282e-01 & 9.346e-01  \\
				&               & (1.35e-01)   & (6.73e-03) & (3.81e-04)&(9.41e-05) \\
				& ${J}_{R_I,G}$ & 2.789e-02    & {\bf2.449e-02}  & 7.407e-01 & 9.602e-01  \\
				&               & (6.52e-03)   & (1.11e-02) & (2.20e-03)& (5.48e-06) \\
				& $L_{R,G}$     & {\bf \color{red} 2.030e-03}     &   &   &\\
				&               & {\bf \color{red} (3.70e-04)}   &   &   &\\
				\Xhline{0.8\arrayrulewidth}
				\multirow{6}{*}{Example~\eqref{ex3}}
				& ${J}_{P_I,G}$ & 1.285e-02    &  {\bf1.027e-03} & 1.910e-00  & 1.472e-00  \\
				&               & (1.60e-02)   & (3.72e-04) & (1.12e-00) &(1.33e-01) \\
				& ${J}_{R_I,G}$ & 1.336e+01    &  {\bf2.397e-03} & 2.435e-00  & 7.683e-01 \\
				&               & (4.97e-00)   & (2.63e-03) & (1.03e-01) &(1.20e-02) \\
				& $L_{R,G}$     & {\bf \color{red} 1.410e-03}     &   &  & \\
				&               & {\bf \color{red} (2.56e-04)}   &   &  & \\
				\Xhline{3\arrayrulewidth}
			\end{tabular}
		\end{center}
	}
	\vskip-.2truecm
	
\end{table}

\paragraph{Activation function and random Fourier feature embedding}
For the activation function, we compare the choices of the $\sin$ and $\tanh$ activation functions. We note again that we used the sine activation function in the previous results. 

Furthermore, we study including random Fourier feature embedding in the fully connected neural network $N({\bf x};\theta)$ as an additional first layer with $2m$ nodes, in addition to the $n$ nodes for the remaining hidden layers, such that
$$
	N_{RFF}({\bf x};\theta):=N({\bf \gamma}({\bf x});\theta),
$$
where
$${\bf \gamma}({\bf x}):=\begin{pmatrix}
	\sin(B {\bf x})\\
	\cos(B {\bf x})
\end{pmatrix}
$$
and each entry of $B \in \mathbb{R}^{m \times 2}$ is sampled from a Gaussian distribution $G(0,\sigma^2 )$ with a user-specified hyper parameter $\sigma$. We also note that $\sigma \in [1,10]$ is recommended in~\cite{wang2023expert}. For the random Fourier feature (RFF) case, we have the resulting network $N_{RFF}({\bf x};\theta)$ with one more layer $\gamma({\bf x})$ of $2m$ output values and it thus has more parameters than that without the RFF, i.e., four hidden layers and $n=35$ nodes per hidden layer.

In~\cref{tab:activation:ex123}, 
we list the results obtained for the test examples~\eqref{ex1} with $k=4$, \eqref{ex2} with $N=4$, and~\eqref{ex3} with $A=100$ and $\varepsilon=0.01$,
with varying the activation function and including of RFF with $\sigma=1$ and $m=17$. 
For the $J_{P,G}$ and $J_{R,G}$ formulations, the smallest error values among the five test cases of $w_B=10^k$, $k=1,\cdots,5$, are reported.
We observe that in Example~\eqref{ex1}, both $\sin$ and $\tanh$ activation functions perform well while the sine activation function gives better results for the examples in~\cref{ex2,ex3} with multi-oscillatory components and high-contrast, oscillatory layers, respectively.

The RFF embedding helps to reduce errors in the examples in~\cref{ex1,ex3}. However, the improvement seems to be problem-dependent, as we could not observe improvements for the second example, Example~\eqref{ex2}. Of course, a variation of the $\sigma$ values may improve the errors, but this simply introduces additional hyper parameter tuning for the $\sigma$ value.

\begin{table}[ht!]
	
	\caption{Activation function and RFF study on Example~\eqref{ex1} with $k=4$, Example~\eqref{ex2} with $N=4$, and Example~\eqref{ex3} with $A=100$ and $\epsilon=0.01$: average and standard deviation of the relative $L^2$-errors (in parenthesis). {The numbers in bold indicate the best result for each loss formulation.}} \label{tab:activation:ex123}
	{\normalsize \renewcommand{\arraystretch}{1.0}
		\begin{center}
			\vskip-.3truecm
			\begin{tabular}{ccccc}
				
				\hline\hline
				&  & Example~\eqref{ex1}  & Example~\eqref{ex2} & Example~\eqref{ex3} \\
				\Xhline{3\arrayrulewidth}
				\multirow{8}{*}{sine}
				& $J_{P,G}$ & {\bf5.100e-04}& {\bf7.620e-03}    & 1.130e-02\\
				&           & (1.66e-04)    &  (3.14e-03)       & (8.97e-03)\\
				& $J_{R,G}$ & {4.917e-02}   & {\bf5.708e-03}    & 2.575e-02\\
				&           & (3.75e-03)    & (4.67e-02)        & (1.51e-03)\\
				& $L_{P,G}$ & {2.420e-03}   & 3.161e-01         & {\bf1.845e-02}\\
				&           & (2.20e-03)    & (4.14e-01)        & (2.00e-03)\\
				& $L_{R,G}$ & {6.400e-04}   & {\bf2.030e-03}    & 1.410e-03\\
				&           & (1.36e-04)    & (3.70e-04)        & (2.56e-04)\\
				\Xhline{0.8\arrayrulewidth}
				\multirow{8}{*}{tanh}
				& $J_{P,G}$ & {1.330e-03}   &  5.787e-02  & 4.057e-02 \\
				&           & (1.44e-04)    &  (2.11e-02) & (1.77e-02)\\
				& $J_{R,G}$ & {3.857e-02}   &  5.537e-01  & 1.00e-00\\
				&           & (1.41e-03)    &  (2.00e-01) & (0.00e-00)\\
				& $L_{P,G}$ & {6.230e-03}   &  3.833e-01  & 3.763e-01\\
				&           & (1.83e-03)    &  (3.27e-01) & (5.17e-01)\\
				& $L_{R,G}$ & {1.320e-03}   &  4.568e-02  & 6.690e-00\\
				&           & (9.80e-05)    &  (2.70e-02) & (3.13e-00)\\
				\Xhline{0.8\arrayrulewidth}
				\multirow{8}{*}{RFF ($\sigma=1$)}
				& $J_{P,G}$ & {7.800e-04}      & 9.450e-03   & {\bf2.090e-03}\\
				&           & (4.82e-04)       & (3.40e-03)  & (4.83e-04)\\
				& $J_{R,G}$ & {\bf1.406e-02}   & 1.147e-01   & {\bf4.430e-03}\\
				&           & (4.36e-03)       & (2.32e-02)  & (1.57e-03)\\
				& $L_{P,G}$ & {\bf1.590e-03}   & {\bf8.124e-02} & 1.943e-02\\
				&           & (4.28e-04)       & (3.63e-02)     & (2.47e-03)\\
				& $L_{R,G}$ & {\bf4.200e-04}   & {6.840e-03}    & {\bf1.030e-03}\\
				&           & (5.59e-05)       & (1.29e-03)     & (2.55e-04)\\
				\Xhline{3\arrayrulewidth}
			\end{tabular}
		\end{center}
	}
	\vskip-.2truecm
	
\end{table}

\begin{table}[ht!]
	
	\caption{Activation function and RFF study on Example~\eqref{ex1} with $k=4$: average and standard deviation of the relative $L^2$-errors (in parenthesis).} \label{tab:activation:ex1}
	{\normalsize \renewcommand{\arraystretch}{1.0}
		\begin{center}
			\vskip-.3truecm
			\begin{tabular}{ccccccc}
				
				\hline\hline
				& $w_B$ & 1  & 10 & 100 & 1000 & 10000 \\
				\Xhline{3\arrayrulewidth}
				\multirow{8}{*}{sine}
				& $J_{P,G}$ & 3.205e-02   & {9.150e-03}    & 1.730e-03    & 7.000e-04    & {\bf5.100e-04} \\
				&           & (6.82e-03) & (1.89e-03) & (5.64e-04) & (5.64e-05) & (1.66e-04) \\
				& $J_{R,G}$ & 1.031e-00   & 2.720e-01 & {\bf4.917e-02}    & 1.000e-00  &1.000e-00 \\
				&           & (4.77e-03) & (1.19e-03) & (3.75e-03) & (1.09e-04) & (3.86e-07) \\
				& $L_{P,G}$ & {\bf2.420e-03}    &   &  &  & \\
				&           & (2.20e-03) &   &  &  & \\
				& $L_{R,G}$ & {\bf6.400e-04}    &   &  &  & \\
				&           & (1.36e-04) &   &  &  & \\
				\Xhline{0.8\arrayrulewidth}
				\multirow{8}{*}{tanh}
				& $J_{P,G}$ & 6.569e-02    & 1.598e-02    & 4.170e-03  & {\bf1.330e-03}   & 1.400e-03 \\
				&           & (1.20e-02)   & (2.35e-03)   & (1.20e-03) & (1.44e-04)       & (1.78e-04) \\
				& $J_{R,G}$ & 5.391e-01    & 2.054e-01    & {\bf3.857e-02}   & 8.064e-01  & 1.000e-00\\
				&           & (3.79e-04)   & (9.81e-04)    & (1.41e-03)      & (3.87e-01) & (0.00e-00) \\
				& $L_{P,G}$ & {\bf6.230e-03}    &   &  &  & \\
				&           & (1.83e-03) &   &  &  & \\
				& $L_{R,G}$ & {\bf1.320e-03}    &  &  &  & \\
				&           & (9.80e-05) &   &  &  & \\
				\Xhline{0.8\arrayrulewidth}
				\multirow{8}{*}{RFF ($\sigma=1$)}
				& $J_{P,G}$ & 3.580e-02    & 6.670e-03    & 1.220e-03    & 5.200e-03 &  {\bf7.800e-04} \\
				&           & (1.26e-02)   & (1.44e-03)   & (1.81e-04)    & (1.63e-04)      & (4.82e-04)\\
				& $J_{R,G}$ & 1.033e-01    & 2.698e-01    & {5.253e-02}   & {\bf1.406e-02}  &  8.036e-00 \\
				&           & (6.20e-03)   & (2.85e-03)    & (4.07e-03)   & (4.36e-03)      & (3.93e-01) \\
				& $L_{P,G}$ & {\bf1.590e-03}    &   &  &  & \\
				&           & (4.28e-04) &   &  &  & \\
				& $L_{R,G}$ & {\bf4.200e-04}    &  &  &  & \\
				&           & (5.59e-05) &   &  &  & \\
				\Xhline{3\arrayrulewidth}
			\end{tabular}
		\end{center}
	}
	\vskip-.2truecm
	
\end{table}

\begin{table}[ht!]
	\caption{Activation function and RFF study on Example~\eqref{ex2} with $N=4$: average and standard deviation of the relative $L^2$-errors (in parenthesis).} \label{tab:activation:ex2}
	{\normalsize \renewcommand{\arraystretch}{1.0}
		\begin{center}
			\vskip-.3truecm
			\begin{tabular}{ccccccc}
				
				\hline\hline
				& $w_B$ & 1  & 10 & 100 & 1000 & 10000 \\
				\Xhline{3\arrayrulewidth}
				\multirow{8}{*}{sine}
				& $J_{P,G}$ & 5.985e-01   & {8.874e-02}    & 4.086e-02    & {\bf7.620e-03}    & {1.878e-02} \\
				&           & (3.41e-01)  & (6.10e-02)     & (2.51e-02)   & (3.14e-03)        & (2.91e-02) \\
				& $J_{R,G}$ & 1.685e-00   &  3.155e-01 & {6.306e-02} & {\bf5.708e-02}  &1.000e-00 \\
				&           & (7.90e-03) & (3.45e-03)  &  (5.61e-03) & (4.67e-02) & (2.75e-06) \\
				& $L_{P,G}$ & {\bf3.161e-01}    &   &  &  & \\
				&           & (4.14e-01) &   &  &  & \\
				& $L_{R,G}$ & {\bf2.030e-03}    &   &  &  & \\
				&           & (3.70e-04) &   &  &  & \\
				\Xhline{0.8\arrayrulewidth}
				\multirow{8}{*}{tanh}
				& $J_{P,G}$ & 1.154e-00    & 8.173e-01    & 2.761e-01     & {1.271e-01}   & {\bf5.787e-02} \\
				&           & (1.89e-01)   & (3.16e-01)   & (4.45e-02)    & (4.80e-02)       & (2.11e-02) \\
				& $J_{R,G}$ & {\bf5.573e-01}    & 8.793e-01    & {8.684e-01}   & 9.883e-01  & 1.000e-00\\
				&           & (2.00e-01)   & (2.39e-01)    & (4.09e-02)   & (2.34e-02) & (2.92e-08) \\
				& $L_{P,G}$ & {\bf3.833e-01}    &   &  &  & \\
				&           & (3.27e-01) &   &  &  & \\
				& $L_{R,G}$ & {\bf4.568e-02}    &  &  &  & \\
				&           & (2.70e-02) &   &  &  & \\
				\Xhline{0.8\arrayrulewidth}
				\multirow{8}{*}{RFF ($\sigma=1$)}
				& $J_{P,G}$ & 8.048e-01    & 2.140e-01    & 6.896e-02    & {2.165e-02}      &  {\bf9.450e-03} \\
				&           & (1.26e-01)   & (4.93e-02)   & (1.08e-02)    & (1.35e-03)      & (3.40e-03)\\
				& $J_{R,G}$ & 1.807e-00    & 3.524e-01    & {\bf 1.147e-01}   & {1.406e-01}     &  8.717e-01 \\
				&           & (1.48e-02)   & (1.73e-02)    & (2.32e-02)   & (2.23e-01)      & (2.10e-01) \\
				& $L_{P,G}$ & {\bf8.124e-02}    &   &  &  & \\
				&           & (3.63e-02) &   &  &  & \\
				& $L_{R,G}$ & {\bf6.840e-03}    &  &  &  & \\
				&           & (1.29e-03) &   &  &  & \\
				\Xhline{3\arrayrulewidth}
			\end{tabular}
		\end{center}
	}
	\vskip-.2truecm
	
\end{table}

\begin{table}[ht!]
	\caption{Activation function and RFF study on Example~\eqref{ex3} with $A=100$ and $\varepsilon=0.01$: average and standard deviation of the relative $L^2$-errors (in parenthesis).} \label{tab:activation:ex3}
	{\normalsize \renewcommand{\arraystretch}{1.0}
		\begin{center}
			\vskip-.3truecm
			\begin{tabular}{ccccccc}
				
				\hline\hline
				& $w_B$ & 1  & 10 & 100 & 1000 & 10000 \\
				\Xhline{3\arrayrulewidth}
				\multirow{8}{*}{$\sin$}
				& $J_{P,G}$ & 2.306e-01   & {9.111e-02}    & 2.303e-02    & {\bf1.130e-02}    & {1.000e-00} \\
				&           & (3.36e-02) & (3.09e-02) & (7.92e-03) & (8.97e-03) & (1.78e-05) \\
				& $J_{R,G}$ & {\bf2.575e-02}   & 4.128e-02 & {4.152e-01}    & 1.000e-00  &1.000e-00 \\
				&           & (1.51e-03)       & (1.30e-02) & (4.78e-01) & (0.00e-00) & (0.00e-00) \\
				& $L_{P,G}$ & {\bf1.845e-02}    &   &  &  & \\
				&           & (2.00e-03) &   &  &  & \\
				& $L_{R,G}$ & {\bf1.410e-03}    &   &  &  & \\
				&           & (2.56e-04) &   &  &  & \\
				\Xhline{0.8\arrayrulewidth}
				\multirow{8}{*}{$\tanh$}
				& $J_{P,G}$ & 1.292e-00    & 2.026e-02    & 3.451e-01  & {\bf4.057e-02}   & 7.744e-01 \\
				&           & (3.42e-01)   & (2.44e-00)   & (5.33e-01) & (1.77e-02)       & (1.35e-00) \\
				& $J_{R,G}$ & 2.480e-00    & 3.611e-00    & {4.409e-00}   & 1.000e-00  & 1.000e-00\\
				&           & (1.15e-00)   & (1.13e-00)    & (2.46e-00)   & (8.92e-08) & (0.00e-00) \\
				& $L_{P,G}$ & {\bf3.763e-01}    &   &  &  & \\
				&           & (5.17e-01) &   &  &  & \\
				& $L_{R,G}$ & {\bf6.690e-00}    &  &  &  & \\
				&           & (3.13e-00) &   &  &  & \\
				\Xhline{0.8\arrayrulewidth}
				\multirow{8}{*}{RFF ($\sigma=1$)}
				& $J_{P,G}$ & 1.500e-01    & 4.102e-02    & 1.370e-02    & {4.890e-03}  &  {\bf2.090e-03} \\
				&           & (2.84e-02)   & (5.04e-03)   & (1.42e-03)    & (7.16e-04)      & (4.83e-04)\\
				& $J_{R,G}$ & 2.849e-02    & 2.857e-02    & {1.756e-02}   & {\bf4.430e-03}  &  8.390e-03 \\
				&           & (6.45e-03)   & (6.32e-03)    & (1.01e-02)   & (1.57e-03)      & (3.51e-03) \\
				& $L_{P,G}$ & {\bf1.943e-02}    &   &  &  & \\
				&           & (2.47e-03) &   &  &  & \\
				& $L_{R,G}$ & {\bf1.030e-03}    &  &  &  & \\
				&           & (2.55e-04) &   &  &  & \\
				\Xhline{3\arrayrulewidth}
			\end{tabular}
		\end{center}
	}
	\vskip-.2truecm
	
\end{table}

\subsection{Study on loss balancing schemes} \label{sec:loss_balancing}

Next, we compare the loss balancing schemes listed in~\cref{tab:loss:balance}
for the test examples in~\cref{ex1,ex2,ex3}.
In our computations, we employ neural networks with $n = 35$ nodes per hidden layer, $n_G = 64$ Gaussian quadrature points, and a total number of $T=100\,000$ epochs.

To improve the training efficiency of PINN, different adaptive weighting methods have been proposed~\cite{mcclenny2023self, maddu2022inverse, wang2023expert}. These methods dynamically adjust the loss weighting factors to balance contributions from different loss components.

\begin{itemize}
	\item \textbf{Self-Adpative Weighting:} 
	The self-adaptive PINN loss function~\cite{mcclenny2023self},  which we denote by the loss function $J_{P,G}({\rm SA})$, introduces an adaptive weight factor that assigns larger weights to training points with higher residual errors. The adaptive weighting method assigns relatively higher importance to regions with larger discrepancies with respect to the underlying PDE or the boundary conditions, which may contribute to improved learning in those areas. Specifically, the loss function is defined as
	$$
		J_{P,G}({\rm SA})(\theta)
		:=
		\int_{\Omega} \lambda_I({\bf x})(\nabla\cdot\nabla U({\bf x};\theta)+f({\bf x}))^2 \, d{\bf x} + \int_{\partial \Omega} \lambda_B({\bf x})(U({\bf x};\theta)-g({\bf x}))^2 \, ds({\bf x}),
	$$
	where $\lambda_k({\bf x}),\ k=I,B$, represent adaptive weight factors which depend on residual magnitudes for the interior and boundary, respectively, and are initially set to $1.0$. We update only the boundary weight factors using the Adam optimizer with a learning rate of $1.0$. {In our experience, updating only the boundary weight factors is more effective for obtaining accurately trained solutions than updating both the interior and boundary weight factors.}
	
	\item \textbf{Inverse-Dirichlet Weighting:} 
	The inverse-Dirichlet weighting method~\cite{maddu2022inverse}, {which we denote by the loss function} $J_{P,G}({\rm invD})$, adjusts the loss weights based on the variance of backpropagated gradients. The standard deviations of the gradients across different loss terms are normalized, which could help {to} balance their contributions during training. The weight update rule is given by
	\begin{align*}
		\hat{w}_k^{(\tau)}
		&=
		\frac{\displaystyle\max_k(\text{std}(\nabla_\theta J_k(\theta^{(\tau)})))}{\text{std}(\nabla_\theta J_k(\theta^{(\tau)}))}, \\
	w_k^{(\tau+1)}&=\alpha w_k^{(\tau)}+(1-\alpha)\hat{w}_k^{(\tau)},
	\end{align*}
	where $k=I,B$ denote the indices corresponding to the interior and boundary loss terms, respectively, and $\tau$ denotes the training epoch. We also set $\alpha=0.5$, as in~\cite{maddu2022inverse}, and the initial value $w_k^{(0)}=1$.
	
	\item \textbf{Gradient-Norm Balancing} 
	The gradient-norm balancing method~\cite{wang2023expert}, {which we denote by the loss function} $J_{P,G}({\rm gradN})$, aims to equalize the norms {of the gradient} of {the different} weighted loss term{s}, effectively ensuring their balance. This approach mitigates the tendency of the model to focus too much on minimizing a specific loss term during training, which could make the optimization process more stable. 
    The weight update rule is given by
	\begin{align*}
		\hat{w}_k^{(\tau)}&=\frac{\lVert{\sum_k\nabla_\theta J_k(\theta^{(\tau)})}\rVert}{\lVert\nabla_\theta J_k(\theta^{(\tau)})\rVert}, \\
		w_k^{(\tau+1)}&=\alpha w_k^{(\tau)}+(1-\alpha)\hat{w}_k^{(\tau)},
	\end{align*}
	where $k=I,B$ denotes the indices corresponding to the interior and boundary loss terms, respectively, and $\tau$ denotes the training epoch. We set $\alpha=0.9$ as in~\cite{wang2023expert}, and the initial value $w_k^{(0)}=1$.
\end{itemize}

For the standard case $J_{P,G}({\rm W})$ and $J_{R,G}({\rm W})$ of a constant user-defined weight $w_B$, we present the smallest error obtained among tests for five $w_B$ values, $w_B=10^k$, $k=0,1,\ldots,4$. In~\cref{tab:loss:bal:result}, the obtained error results are listed for the three test examples. For the Example~\eqref{ex1}, the PINN loss $J_{P,G}({\rm W})$ with a large weight factor and the deep Ritz loss with the augmented Lagrangian term, $L_{R,G}$, perform well, while for the other two Examples~\eqref{ex2} and~\eqref{ex3}, the case of $L_{R,G}$ performs the best among the proposed loss balancing schemes. Among the many weighting schemes tested, none performed better overall than the constant weighting scheme for PINNs, $J_{P,G}({\rm W})$, and the augmented Lagrangian approach for the deep Ritz method, $L_{R,G}$.

\begin{table}
	\begin{center}
		\begin{tabular}{|l|l|}
			\hline
			Notation & Loss balancing schemes \\
			\hline
			$J_{P,G}({\rm W})$ & PINN loss with (large) constant weight factor $w_B$ \\			$J_{P,G}({\rm SA})$ & PINN loss with a self-adaptive weight factor \\
			$J_{P,G}({\rm invD})$ & PINN loss with an inverse-Dirichlet weight \\
			$J_{P,G}({\rm gradN})$ & PINN loss with {a gradient norm}\\
			$L_{P,G}$ & PINN loss with an augmented Lagrangian \\
			$J_{R,G}({\rm W})$ & Deep Ritz loss with (large) constant weight factor $w_B$ \\
			$L_{R,G}$ & Deep Ritz loss with an augmented Lagrangian \\
			\hline
		\end{tabular}
		\caption{Notation for the loss balancing scheme study.
		}\label{tab:loss:balance}
	\end{center}
\end{table}

\begin{table}[ht!]
	\caption{Loss balancing scheme study on test examples in \cref{ex1,ex2,ex3}: average and standard deviation of the relative $L^2$-errors (in parenthesis).} \label{tab:loss:bal:result}
	{\normalsize \renewcommand{\arraystretch}{1.0}
		\begin{center}
			\vskip-.3truecm
			\begin{tabular}{ccccccc}
				
				\hline\hline
				Schemes & Rank & Example~\eqref{ex1} & Rank & Example~\eqref{ex2}& Rank & Example~\eqref{ex3} \\
				\Xhline{3\arrayrulewidth}
				$J_{P,G}({\rm W})$   &  1  & {\bf5.100e-04}   & 2  & 1.260e-02  & 5 & 1.605e-02   \\
				&    & (1.66e-04)   &   & (4.34e-03)  &  & (2.71e-03)  \\
				$J_{P,G}({\rm SA})$  &   5 & 1.220e-03   & 3  & 5.093e-02  & 4 & 8.281e-03   \\
				&    & (2.73e-04)   &   & (5.65e-03)  &  & (1.30e-03)  \\
				$J_{P,G}({\rm invD})$&   4 & 8.800e-04   & 7  & 3.681e-01  & 3 & 7.370e-03   \\
				&    & (3.00e-04)   &   & (2.59e-01)  &  & (1.19e-03)  \\
				$J_{P,G}({\rm gradN})$ & 3 &8.400e-04    & 4  & 6.369e-02  & 2 & 7.010e-03  \\
				&  & (2.26e-04)   &   & (4.09e-02)  &  & (1.36e-03) \\
				$L_{P,G}$        & 6 & 2.420e-03   & 6  & 1.345e-01  & 6 & 2.615e-02   \\
				&  & (2.20e-03)   &   & (2.99e-02)  &  & (2.07e-03) \\
				$J_{R,G}({\rm W})$     & 7 & 4.917e-02   & 5  & 1.028e-01  & 7 & 5.511e-02   \\
				&  & (3.95e-03)   &   & (2.37e-02)  &  & (1.02e-03) \\
				$L_{R,G}$        & 2 & 6.400e-04   & 1  & {\bf6.610e-03}  & 1 & {\bf4.170e-03} \\
				&  & (1.36e-04)   &   & (2.28e-03)  &  & (5.28e-04)\\
				\Xhline{3\arrayrulewidth}
			\end{tabular}
		\end{center}
	}
	\vskip-.2truecm
	
\end{table}

\subsection{Study on optimizers} \label{sec:optimizers}

Finally, we compare the performance of the four loss formulations, $J_{P,G}$, $J_{R,G}$, $L_{P,G}$, and $L_{R,G}$ depending on the optimizer choice. In particular, we consider the following three settings:
\begin{itemize}
	\item Adam~\cite{adam}
	\item Limited-memory Broyden--Fletcher--Goldfarb--Shanno (L-BFGS)~\cite{liu_limited_1989}
	\item Adam+L-BFGS
\end{itemize}
In the Adam and the L-BFGS cases, we train the network parameters for $T=100\,000$ epochs with a learning rate of $\epsilon=0.001$. For the Adam+L-BFGS case, we first train using the Adam optimizer up to $T_A=80\,000$ epochs and then switch to the L-BFGS optimizer until the final epoch, $T=100\,000$.

The error results obtained are listed in~\cref{tab:optimizer:ex123}, where we, again, consider the test examples in \cref{ex1,ex2,ex3} with $k=1$, $N=4$, and $A=100$ and $\varepsilon=0.01$, respectively. For all examples, we choose a network width of $n=35$ nodes per hidden layer and $n_G=64$ Gaussian quadrature points in each direction. For all the test examples, we can observe that the training process with Adam is stable and the obtained results are more accurate compared to the other optimizer settings.
Furthermore, we note that, for Example~\eqref{ex1}, we obtain good results for all the optimizer choices except for the combination of the $J_{R,G}$ loss formulation and the L-BFGS optimizer. 
Moreover, the $L_{R,G}$ loss yields reasonable convergence for all the optimizer choices.

The loss formulation $L_{R,G}$ seems more robust to both the test examples
and the optimizers compared to other loss formulations.

\begin{table}[ht!]
	\caption{Optimizer study on examples in \cref{ex1}-\eqref{ex3} with $k=1$, $N=4$, and $A=100$ and $\varepsilon=0.01$, respectively: the average of the relative $L^2$-errors, and the standard deviation (in parenthesis), excluding non-converging seeds, the symbol - indicates that the optimizer does not give a convergent solution.} \label{tab:optimizer:ex123}
	{\normalsize \renewcommand{\arraystretch}{1.0}
		\begin{center}
			\vskip-.3truecm
			\begin{tabular}{ccccc}
				
				\hline\hline
				&  & Adam  & L-BFGS & Adam+L-BFGS \\
				\Xhline{3\arrayrulewidth}
				\multirow{8}{*}{Example~\eqref{ex1}}
				& $J_{P,G}$ & 3.000e-04   & {4.570e-05}  &   5.620e-05  \\
				&           & (4.85e-05)  & (0.00e-00)   &   (2.97e-05) \\
				& $J_{R,G}$ & 7.060e-03   & 2.302e-01    &   6.530e-03  \\
				&           & (9.74e-04)  & (2.70e-06)   &   (1.19e-03) \\
				& $L_{P,G}$ & {1.500e-04} & 1.225e-04    &   5.069e-04  \\
				&           & (4.94e-05)  & (8.32e-05)   &   (2.43e-04) \\
				& $L_{R,G}$ & {4.900e-04} & 9.230e-05    &   5.168e-04  \\
				&           & (2.35e-04)  & (2.22e-05)   &   (1.98e-04) \\
				\Xhline{0.8\arrayrulewidth}
				\multirow{8}{*}{Example~\eqref{ex2}}
				& $J_{P,G}$ & 1.260e-02    & 4.262e-01   &      -       \\
				&           & (4.34e-03)   & (8.88e-02)  &      -       \\
				& $J_{R,G}$ & 1.028e-01    & 3.401e-01   &   4.653e-01  \\
				&           & (2.37e-02)   & (4.95e-04)  &   (1.80e-01) \\
				& $L_{P,G}$ & 1.345e-01    & 9.161e-00   &      -       \\
				&           & (2.99e-02)   & (1.55e-00)  &      -       \\
				& $L_{R,G}$ & 6.610e-03    & 1.238e-02   &   9.245e-03  \\
				&           & (2.28e-03)   & (1.31e-02)  &   (3.03e-03) \\
				\Xhline{0.8\arrayrulewidth}
				\multirow{8}{*}{Example~\eqref{ex3}}
				& $J_{P,G}$ & 1.130e-02    & -           &   4.825e-03  \\
				&           & (8.97e-03)   & -           &   (0.00e-00) \\
				& $J_{R,G}$ & 2.575e-02    & 1.711e-02   &   2.438e-02  \\
				&           & (1.51e-03)   & (7.60e-05)  &   (7.45e-04) \\
				& $L_{P,G}$ & 1.845e-02    & -           &      -       \\
				&           & (2.00e-03)   & -           &      -       \\
				& $L_{R,G}$ & 1.410e-03    & 4.772e-04   &   1.512e-03  \\
				&           & (2.56e-04)   & (2.65e-05)  &   (2.33e-04) \\
				\Xhline{3\arrayrulewidth}
			\end{tabular}
		\end{center}
	}
	\vskip-.2truecm
	
\end{table}

\section{Some additional challenging examples}\label{sec:DeepRitz:numerics}

In~\cref{numerics:result}, we have focused on variations of Poisson model problems. Finally, in this section, we consider some challenging examples:
\begin{itemize}
	\item{\bf three-dimensional problems}
	\item{\bf nonlinear $p$-Laplacian problem with increasing $p$ }
	\item{\bf eigenvalue problem}
\end{itemize}
We will observe that deep Ritz formulation with augmented Lagrangian loss function, i.e., $L_{R,G}$, which had already performed well in the results reported in the previous section, outperforms the other approaches considered for these more challenging cases.

\subsection{Three-dimensional problems}

In this subsection, we consider three-dimensional Poisson model problems, comparing the $J_{P,G}$, $J_{R,G}$, $L_{P,G}$, and $L_{R,G}$ approaches in terms of computing times and solution accuracy. In particular, we consider the following three solutions
\begin{eqnarray}
	u(x,y,z) & = & \sin(k \pi x)\sin(k \pi y) \sin(k \pi z), \label{3d:ex1} \\
	u(x,y,z) & = & \frac{1}{N}\sum_{\ell=1}^N \sin(2^\ell \pi x) \sin(2^\ell \pi y) \sin(2^\ell \pi z), \label{3d:ex2} \\
	u(x,y,z) & = & Ax(1-x)y(1-y)z(1-z) \sin(\frac{(x-0.5)(y-0.5)(z-0.5)}{\varepsilon}), \label{3d:ex3}
\end{eqnarray}
for $(x,y,z) \in \Omega$. We then choose the right hand side $f$ and boundary function $g$ in~\cref{model:poisson} accordingly. The values $k$, $N$, and $A$, $\varepsilon$ are chosen as
$$
	k=4,\; N=2,\; A=100, \varepsilon=0.01,
$$
and the domain $\Omega$ is a unit cubic domain, i.e., $\Omega=(0 \, 1)^3$. 

For all the test examples, we consider a neural network with $n=100$ nodes per hidden layer and a sampling set with $n_G=32$ Gaussian quadrature points in each direction. We train the network parameters for $T=100\,000$ training epochs and report error values based on the minimum error indicator throughout the whole training process.

In~\cref{tab:3d:result}, the average and standard deviation of the relative $L^2$-error values are listed for the PINN and deep Ritz formulations. 
The results are obtained from five different parameter initializations. For the cases, $J_{P,G}$ and $J_{R,G}$, the weight factor $w_B$ are set to $w_B=10^{k}$ with $k=0,1,\ldots,4$ and the minimum error values are reported among the five different $w_B$ cases. For the cases, $L_{P,G}$ and $L_{R,G}$, the weight factor $w_B$ is simply set to 1 since the augmented Lagrangian term is included to deal with the imbalance between the differential equation and the boundary condition terms in the loss function. For the test examples~\eqref{3d:ex1} and \eqref{3d:ex2}, the PINN formulation $J_{P,G}$ with a large weight factor gives the smallest error values but with a much more computation time than in the deep Ritz formulations, $J_{R,G}$ and $L_{R,G}$. The deep Ritz formulation $L_{R,G}$ gives comparable error results to those obtained from $J_{P,G}$ with about two or three factors larger error values. For the test example~\eqref{3d:ex3}, the deep Ritz formulation $L_{R,G}$ gives the smallest error values with a much lesser computation time than in the PINN formulation $J_{P,G}$ with a larger weight factor.
In addition, the advantage in the deep Ritz formulation $L_{R,G}$ is no additional tuning for the hyper parameter $w_B$, while the performance of $J_{P,G}$ highly depends on the choice of $w_B$.

\begin{table}
	\begin{center}
		\begin{tabular}{|l|r|r|r|r|}
			\hline
			& Example~\eqref{3d:ex1} & Example~\eqref{3d:ex2} & Example~\eqref{3d:ex3} & Computation time\\
			\hline
			$J_{P,G}$ &  {\bf 2.440e-03}   &   {\bf 1.950e-03}  & 7.840e-03  & 2\,300\,s\\
			& (4.43e-04)   &   (2.56e-04) & (1.95e-03) &  \\
			$J_{R,G}$ &  2.537e-01   &   8.080e-02  & 1.408e-01  & {\bf 620s} \\
			& (8.56e-03)   &   (1.02e-02) & (1.08e-02) &  \\
			$L_{P,G}$ &  3.179e-02   &   5.050e-03  & 1.445e-02  & 2\,300\,s \\
			& (1.21e-02)   &   (3.19e-03) & (2.20e-04) &  \\
			$L_{R,G}$ & 8.790e-03    &   4.870e-03  & {\bf 4.760e-03}  & {\bf 620\,s} \\
			& (2.28e-03)   &   (1.14e-03) & (1.01e-03) &  \\
			\hline
		\end{tabular}
	\end{center}
	\caption{Error and computation time results for three-dimensional test examples in~\cref{3d:ex1,3d:ex2,3d:ex3}: average and standard deviation of the relative $L^2$-errors and average computing times; best results in boldface.}
	\label{tab:3d:result}
\end{table}

\subsection{$p$-Laplacian problem}

Next, we consider a $p$-Laplacian problem with a smooth solution,
\begin{equation} \label{eq:p lap}
	\begin{split}
		-\Delta_p u&=f\quad \text{in } \Omega:=(0,1)^2, \\
		u&=g\quad \text{on}\, \partial\Omega,
	\end{split}
\end{equation}
where $f$ and $g$ are chosen such that the exact solution is
\begin{equation*}\label{ex:p-Lap}
	u^*(x,y)=\sin(2\pi x)\sin(2\pi y).
\end{equation*}
Here, the $p$-Laplace operator is defined as $\Delta_{p} u := {\rm div} (|\nabla u |^{p-2} \nabla u)$; note that the $p$-Laplacian simplifies to the standard Laplacian, which we considered in the previous model problems, for the case $p=2$.

For the $p$-Laplacian model problem, the PINN formulation of the loss function reads
\begin{equation*}\label{loss:PNG_p_lap}
	\begin{split}
		J_{P,G}(\theta)&:=w_I \sum_{ {\bf x} \in X_G(\Omega)} \left(\nabla \cdot \left(\lvert\nabla U({\bf x};\theta)\rvert^{p-2}\nabla U({\bf x};\theta)\right)+f({\bf x})\right)^2 w({\bf x}) \\
		&+w_B \sum_{ {\bf x} \in X_G(\partial \Omega)} \left(U({\bf x};\theta)-g({\bf x}) \right)^2 w({\bf x}),
	\end{split}
\end{equation*}
and the deep Ritz formulation reads
\begin{equation*}\label{loss:DRG_p_lap}
	\begin{split}
		J_{R,G}(\theta)&:=w_I \sum_{ {\bf x} \in X_G(\Omega)} \left( \frac{1}{p}|\nabla U({\bf x};\theta)|^p-f({\bf x})U({\bf x};\theta) \right) w({\bf x})\\
		&+w_B \sum_{ {\bf x} \in X_G(\partial \Omega)} (U({\bf x};\theta)-g({\bf x}))^2 w({\bf x}).
	\end{split}
\end{equation*}

For our numerical experiments, we employ a fully connected neural network with width $n=35$, Gaussian quadrature with $n_G=64$ sampling points in each direction, and a total number of $T=100\,000$ training epochs. In the cases of $J_{P,G}$ and $J_{R,G}$, we tested five difference choices for the weight factor $w_B=10^k$ with $k=0,1,\ldots,4$, and we report the minimum $L^2$-error among the five cases. To deal with the increasing magnitude of $f$ for the higher values of $p$, we adjust the weight factor $w_I$: 
\begin{equation}\label{wI:p-L:weight}
	w_I=\frac{1}{{\int_{\Omega}\lvert f(x,y)\rvert \, d{\bf x}}}.
\end{equation}

In~\cref{tab:p_lap_34}, we report the relative $L^2$-errors for the $p$-Laplacian model problem in~\cref{eq:p lap} with increasing $p$ values, $p=3,4,\ldots,7$, and the weight factor $w_I$ defined in~\cref{wI:p-L:weight}. For $p=3$, the PINN formulations, $J_{P,G}$ and $L_{P,G}$, give more accurate results than the respective deep Ritz formulations, $J_{R,G}$ and $L_{R,G}$. For larger $p$ values, the deep Ritz formulations $J_{R,G}$ and $L_{R,G}$ yield smaller errors. Moreover, the deep Ritz methods appear to be more robust towards increasing values of $p$, in the sense that the error increase is less strong. In both PINN and deep Ritz formulations, the the augmented Lagrangian formulation helps to reduce the errors.

To show the importance of the hyper parameter choice $w_I$, we also present the error results for $w_I=1$ in~\cref{tab:p_lap_34:wI:one}. Only for $p=3$, the simple choice $w_I=1$ yields smaller errors in the $J_{P,G}$ and $L_{R,G}$ cases, while the error results are worse than those for~\cref{wI:p-L:weight}, as listed in~\cref{tab:p_lap_34}.

\begin{table}[ht!]
	\caption{$p$-Laplacian problem in~\eqref{eq:p lap} for increasing $p$ and with $w_I$ as defined in~\cref{wI:p-L:weight}: the average and standard deviation (in parenthesis) of the relative $L^2$-errors of the four methods. The best result for each value of $p$ is in boldface.}\label{tab:p_lap_34}
	{\normalsize \renewcommand{\arraystretch}{1.0}
		\begin{center}
			\vskip-.3truecm
			\begin{tabular}{cccccc}
				
				\hline\hline
				$p$      & 3  	 & 	 4    & 5 & 6 & 7 \\
				\Xhline{3\arrayrulewidth}
				$J_{P,G}$   &   7.991e-04    &   1.000e-00    &   1.000e-00    &   1.000e-00   &   1.000e-00     \\
				&   (1.91e-04)   &   (3.73e-07)   &   (2.56e-06)   &   (4.70e-06)   &   (5.97e-04)    \\
				$J_{R,G}$   &   2.048e-02    &   2.396e-02    &   2.242e-02    &   2.041e-02    &   2.560e-02     \\
				&   (3.29e-03)   &   (6.72e-03)   &   (3.83e-03)   &   (3.61e-03)   &   (4.68e-03)    \\
				$L_{P,G}$   &  {\bf5.549e-04}& {\bf8.791e-04}  &   1.704e-02    &   1.350e-01    &   1.921e-01     \\
				&   (3.26e-04)   &   (1.68e-04)   &   (2.46e-02)   &   (2.18e-01)   &   (3.78e-01)    \\
				$L_{R,G}$   &  1.834e-03 & {2.102e-03} & {\bf2.471e-03} & {\bf2.669e-03} & {\bf2.736e-03}  \\
				&   (5.11e-04)   &   (4.12e-04)   &   (1.11e-03)   &   (7.65e-04)   &   (1.17e-03)    \\
				\Xhline{3\arrayrulewidth}
			\end{tabular}
		\end{center}
	}
	\vskip-.2truecm
\end{table}

\begin{table}[ht!]
	\caption{$p$-Laplacian problem in~\eqref{eq:p lap} with increasing $p$ and $w_I=1$: the average and standard deviation (in parenthesis) of the relative $L^2$-errors of the four methods. The best result for each value of $p$ is in boldface.}\label{tab:p_lap_34:wI:one}
	{\normalsize \renewcommand{\arraystretch}{1.0}
		\begin{center}
			\vskip-.3truecm
			\begin{tabular}{cccccc}
				
				\hline\hline
				$p$      & 3  	 & 	 4    & 5 & 6 & 7 \\
				\Xhline{3\arrayrulewidth}
				$J_{P,G}$   &   7.794e-04    &   1.000e-00    &   7.044e-01    &   1.000e-00   &   1.021e-00     \\
				&   (3.38e-04)   &   (7.71e-07)   &   (5.74e-01)   &   (2.36e-06)   &   (2.54e-02)    \\
				$J_{R,G}$   &   2.891e-02    &   2.392e-02    &  {\bf6.443e-02}&  {\bf2.792e-01}& {\bf8.593e-01} \\
				&   (2.74e-03)   &   (7.19e-03)   &   (2.54e-03)   &   (2.45e-03)   &   (3.14e-03)    \\
				$L_{P,G}$   &  {1.012e-01}& {4.216e-00}  &   2.119e-02    &   1.971e-00    &   1.627e-00     \\
				&   (6.74e-04)   &   (1.34e-00)   &   (1.67e-00)   &   (1.58e-00)   &   (1.10e-00)    \\
				$L_{R,G}$   &  {\bf3.458e-04} & {\bf6.507e-03} & {9.227e-02} & {1.361e-00} & {1.464e-00}  \\
				&   (5.36e-05)   &   (2.00e-03)   &   (2.28e-03)   &   (1.16e-02)   &   (1.18e-02)    \\
				\Xhline{3\arrayrulewidth}
			\end{tabular}
		\end{center}
	}
	\vskip-.2truecm
\end{table}

\subsection{Eigenvalue problem}

In this subsection, we consider the eigenvalue problem
\begin{equation} \label{eq:eigen}
	\begin{split}
		-\Delta u +vu&=  \mu  u \quad \text{in}\, \Omega, \\
		u&=0\quad \text{on}\, \partial\Omega,
	\end{split}
\end{equation}
{where $v$ is a given potential function and $\mu$ is an eigenvalue.} It is well-known that the smallest eigenvalue $\lambda_{\min}$ minimizes the following functional, called the Rayleigh quotient:
\begin{equation*}
	\lambda_{\min}
	=
	\min_{u\rvert_{\partial \Omega}=0} \quad \frac{\int_\Omega \lvert\nabla u\rvert^2 dx + \int_\Omega vu^2 dx}{\int_\Omega u^2 dx}.
\end{equation*}
To avoid the case of the trivial solution, i.e., $u\equiv0$, we form the following constrained minimization problem with an additional constraint, $\int_{\Omega} u^2 d {\bf x}=1$:
\begin{equation*}
	\begin{split}
		\min_{\substack{u\rvert_{\partial \Omega}=0, \\ \int_\Omega u^2 d {\bf x} = 1}} &\quad \frac{\int_\Omega \lvert\nabla u\rvert^2 dx + \int_\Omega vu^2 dx}{\int_\Omega u^2 dx}.
	\end{split}
\end{equation*}

In our computations, we use the following loss function, augmenting the constraints
$u|_{\partial \Omega}=0$ and $\int_{\Omega} u^2 d {\bf x}=1$  with Lagrange multipliers $\lambda({\bf x})$ and $\lambda_C$, respectively:
\begin{equation*}
	\begin{split}
		L_{R}(\theta,\lambda,\lambda_C) &= \frac{\int_\Omega \lvert\nabla U({\bf x};\theta)\rvert^2\,d{\bf x} + \int_\Omega v\left(U({\bf x};\theta)\right)^2\,d{\bf x}}{\int_\Omega \left(U({\bf x};\theta)\right)^2\,d{\bf x}}
		+ w_B \int_{\partial \Omega}(U({\bf x};\theta))^2\,ds({\bf x}) \\
		&+ \int_{\partial \Omega}\lambda({\bf x}) U({\bf x};\theta)\,ds({\bf x}) + w_C\left(\int_\Omega \left(U({\bf x};\theta)\right)^2\,d{\bf x} - 1\right)^2\\
		&+ \lambda_C\left(\int_\Omega \left(U({\bf x};\theta)\right)^2\,d{\bf x} - 1\right),
	\end{split}
\end{equation*}
where $w_B$ and $w_C$ are the weight factors associated with the two constraint conditions. As before, we employ a fully connected neural network $U({\bf x};\theta)$ with width $n=35$ nodes per hidden layer, and the loss $L_R(\theta,\lambda,\lambda_C)$ is approximated by $L_{R,G}(\theta,\lambda,\lambda_C)$ using the Gaussian quadrature with $n_G=64$ in each direction. For the Lagrange multipliers $\lambda_C$ and $\lambda({\bf x})$, we set the initial value as $1.0$. We train the network parameters $\theta$ and the Lagrange multipliers $\lambda({\bf x})$ and $\lambda_C$ for a total of $T=100\,000$ epochs with the same learning rates as before.
Since the approximate solution oscillates
over the training epochs, we compute the average value over the last $10\,000$ training epochs to give a stable approximate solution. 


In the following, we consider the two test examples in~\cite{yu2018deep}.

\paragraph{Infinite potential well} We can reformulate the eigenvalue problem in~\cref{eq:eigen} into the following equivalent problem:
\begin{equation} \label{eq:eigen_1}
	\begin{split}
		-\Delta u({\bf x}) &= \mu u({\bf x}), \quad {\bf x}\in \Omega :=(0,1)^2,\\
		u({\bf x})&=0, \quad\qquad {\bf x}\in \partial \Omega,
	\end{split}
\end{equation}
where the potential function is given as
\begin{equation*}
	v(x)=
	\begin{cases}
		0, & x \in [0,1]^2, \\
		\infty, & x\notin [0,1]^2.
	\end{cases}
\end{equation*}
In the case, the smallest nonzero eigenvalue is $\mu_0=2\pi^2$.

In~\cref{tab:eigen_1}, we report the average and standard deviation of the relative errors of the approximate eigenvalue. The training results are affected by both weights $w_B$ and $w_C$. For $w_B$ in the range between $10$ and $100$ and $w_C$ in the range between $10$ and $1\,000$, we obtained the average error values less than $10^{-4}$.

\paragraph{The harmonic oscillator}

Finally, we consider the eigenvalue problem in~\cref{eq:eigen} with the potential function $v({\bf x})=\lvert {\bf x} \rvert^2$,
\begin{equation} \label{eq:eigen_2}
	\begin{split}
		-\Delta u({\bf x}) + \lvert {\bf x} \rvert^2 u &= \mu u({\bf x}), \quad {\bf x}\in \Omega :=(-3,3)^2,\\
		u({\bf x})&=0, \quad\qquad {\bf x}\in \partial \Omega.
	\end{split}
\end{equation}
In this case, the smallest nonzero eigenvalue is $\mu_0=2$.

In~\cref{tab:eigen_2}, the average and standard deviation of the relative errors to the smallest eigenvalue for the eigenvalue problem in~\cref{eq:eigen_2} are reported. Similarly to the previous example, we can observe that the accuracy of the trained results depends on the weight factors $w_B$ and $w_C$. For both $w_B$ and $w_C$ in the range between $10$ and $100$, the average errors are below $10^{-3}$.

In summary, as reported in~\cref{tab:eigen_1,tab:eigen_2}, we obtained very accurate approximate values of the minimum nonzero eigenvalue with the choice of weight factors $w_B$ and $w_C$ in a relatively mild range between $10$ and $100$.

\begin{table}[ht!]
	\caption{Eigenvalue problem in~\eqref{eq:eigen_1}: The average and standard deviation (in parenthesis) of the relative errors depending on the choice of $w_B$ and $w_C$. Relative errors below $10^{-4}$ are marked in boldface.}\label{tab:eigen_1}
	{\normalsize \renewcommand{\arraystretch}{1.0}
		\begin{center}
			\vskip-.3truecm
			\begin{tabular}{@{}cccccc}
				
				\hline\hline
				\backslashbox{$w_B$}{$w_C$} & 1  & 10 & 100 & 1000 & 10000 \\
				\Xhline{3\arrayrulewidth}
				1     & 9.668e-01 & 7.087e-01 & 7.074e-01 & 6.178e-01 & 5.069e-01 \\
				& (7.29e-03)& (1.48e-02)& (1.24e-02)& (8.56e-03)& (1.92e-02)\\
				10    & 5.616e-01 & {\bf1.121e-05} & {\bf1.266e-05} & {\bf3.842e-05} & 1.838e-04 \\
				& (2.89e-03)& (2.61e-06)& (2.93e-06)& (2.07e-05)& (7.32e-05)\\					
				100   & 3.403e-01 & {\bf1.353e-05} & {\bf1.840e-05} & {\bf5.519e-05} & 2.680e-04 \\
				& (8.15e-02)& (1.08e-05)& (8.99e-06)& (4.67e-05)& (2.05e-04)\\
				1000  & 6.137e-02 & 2.476e-03 & 2.078e-03 & 2.629e-03 & 1.262e-03 \\
				& (4.53e-02)& (1.25e-03)& (3.87e-04)& (5.76e-05)& (4.38e-04)\\
				10000 & 8.557e-02 & 1.071e-02 & 5.475e-03 & 2.687e-03 & 1.368e-03 \\
				& (4.39e-02)& (7.76e-03)& (1.72e-03)& (1.85e-03)& (6.42e-04)\\
				\Xhline{3\arrayrulewidth}
			\end{tabular}
		\end{center}
	}
	\vskip-.2truecm
\end{table}

\begin{table}[ht!]
	\caption{Eigenvalue problem in~\eqref{eq:eigen_2}: The average and standard deviation (in parenthesis) of the relative errors depending on the choice of $w_B$ and $w_C$. Relative errors below $10^{-3}$ are marked in boldface.}\label{tab:eigen_2}
	{\normalsize \renewcommand{\arraystretch}{1.0}
		\begin{center}
			\vskip-.3truecm
			\begin{tabular}{cccccc}
				
				\hline\hline
				\backslashbox{$w_B$}{$w_C$} & 1  & 10 & 100 & 1000 & 10000 \\
				\Xhline{3\arrayrulewidth}
				1     & 3.249e-00 & {\bf8.500e-04} & 1.502e-02 & 3.074e-01 & 8.581e-01 \\
				& (2.54e-01)& (1.43e-05)& (1.64e-02)& (2.21e-01)& (4.60e-01)\\
				10    & 2.596e-00 & {\bf8.472e-04} & {\bf9.144e-04} & 1.210e-01 & 1.945e-01 \\
				& (2.15e-00)& (1.94e-05)& (1.15e-05)& (8.33e-02)& (5.81e-02)\\					
				100   & 1.335e-00 & {\bf8.751e-04} & {\bf9.203e-04} & 9.196e-02 & 3.177e-01 \\
				& (1.26e-00)& (2.43e-05)& (3.46e-05)& (7.85e-02)& (2.27e-01)\\
				1000  & 3.667e-00 & 1.231e-03 & 1.074e-03 & 1.442e-03 & 3.324e-01 \\
				& (3.43e-01)& (1.24e-04)& (5.78e-05)& (1.35e-04)& (1.97e-01)\\
				10000 & 8.945e-01 & 3.045e-03 & 2.836e-03 & 2.150e-03 & 4.590e-03 \\
				& (2.39e-01)& (5.33e-04)& (3.05e-04)& (2.86e-04)& (1.18e-03)\\					
				\Xhline{3\arrayrulewidth}
			\end{tabular}
		\end{center}
	}
	\vskip-.2truecm
\end{table}

\section{Conclusions}\label{sec:Conclude}

In this work, extensive numerical studies on hyper parameter choices in neural network approximation of partial differential equations were conducted. While generally applicable rules are out-of-reach to derive, we aim at making some practical suggestions for hyper parameter choices for test examples with typical properties and varying complexity, i.e., smooth solution, multi-component oscillatory solution, and high-contrast, oscillatory interior layer solution. We consider the two most popular formulations of PDE loss functions, the PINN and deep Ritz methods, and compared them for various hyper parameter settings. We have observed that the use of the augmented Lagrangian approach for balancing the PDE and boundary loss terms as well as more accurate numerical integration schemes can improve the performance of deep Ritz formulation. We have observed that, using those techniques, the deep Ritz methods appears to be more accurate and stable compared to the PINN method for the more complex model problems, i.e., the multi-component oscillatory and the high-contrast, oscillatory interior layer cases. The study on various loss balancing schemes indicated good performance of the augmented Lagrangian approach, also being more robust to the model problem complexity. Finally, we have observed that the deep Ritz formulation with the augmented Lagrangian term and with a more accurate integration scheme generally outperforms the other approaches for even more challenging examples, including three-dimensional, nonlinear, and eigenvalue problems, in terms of accuracy and computing time.

Based on our hyper parameter study, our overall suggestion is the following: when a more accurate numerical integration scheme like a Gaussian quadrature is {applicable}
and the model problem can be reformulated as a minimization problem, the deep Ritz formulation with the augmented Lagrangian term and with the quadrature sampling points yield good setup for the hyper parameters. Otherwise, the PINN method is more flexible, in the sense that it can be applied to a wider range of model problems, and its performance seems to be less sensitive to the hyper parameter settings.

\bibliographystyle{plain}
\bibliography{ddforpinn}

\begin{thebibliography}{10}

\bibitem{abadi_tensorflow_2016}
M.~Abadi, A.~Agarwal, P.~Barham, E.~Brevdo, Z.~Chen, C.~Citro, G.~S. Corrado,
  A.~Davis, J.~Dean, M.~Devin, S.~Ghemawat, I.~Goodfellow, A.~Harp, G.~Irving,
  M.~Isard, Y.~Jia, R.~Jozefowicz, L.~Kaiser, M.~Kudlur, J.~Levenberg, D.~Mane,
  R.~Monga, S.~Moore, D.~Murray, C.~Olah, M.~Schuster, J.~Shlens, B.~Steiner,
  I.~Sutskever, K.~Talwar, P.~Tucker, V.~Vanhoucke, V.~Vasudevan, F.~Viegas,
  O.~Vinyals, P.~Warden, M.~Wattenberg, M.~Wicke, Y.~Yu, and X.~Zheng.
\newblock {TensorFlow}: {Large}-{Scale} {Machine} {Learning} on {Heterogeneous}
  {Distributed} {Systems}, March 2016.
\newblock arXiv:1603.04467 [cs].

\bibitem{basir2023investigating}
S.~Basir.
\newblock Investigating and mitigating failure modes in physics-informed neural
  networks (pinns).
\newblock {\em Commun. Comput. Phys.}, 33(5):1240--1269, 2023.

\bibitem{berg2018}
J.~Berg and K.~Nystr{\"o}m.
\newblock A unified deep artificial neural network approach to partial
  differential equations in complex geometries.
\newblock {\em Neurocomputing}, 317:28--41, 2018.

\bibitem{blechschmidt_three_2021}
J.~Blechschmidt and O.~G. Ernst.
\newblock Three ways to solve partial differential equations with neural
  networks — {A} review.
\newblock {\em GAMM-Mitt.}, 44(2):e202100006, 2021.

\bibitem{bradbury_jax_2018}
J.~Bradbury, R.~Frostig, P.~Hawkins, M.~J. Johnson, C.~Leary, D.~Maclaurin,
  G.~Necula, A.~Paszke, J.~VanderPlas, S.~Wanderman-Milne, and Q.~Zhang.
\newblock {JAX}: composable transformations of {Python}+{NumPy} programs, 2018.

\bibitem{cai_physics-informed_2021}
S.~Cai, Z.~Mao, Z.~Wang, M.~Yin, and G.~E. Karniadakis.
\newblock Physics-informed neural networks ({PINNs}) for fluid mechanics: a
  review.
\newblock {\em Acta Mech. Sin.}, 37(12):1727--1738, 2021.

\bibitem{chae2023two}
J.~Chae, K.~Kim, and D.~Kim.
\newblock Two-timescale extragradient for finding local minimax points.
\newblock {\em arXiv preprint arXiv:2305.16242}, 2023.

\bibitem{cuomo_scientific_2022}
S.~Cuomo, V.~S. Di~Cola, F.~Giampaolo, G.~Rozza, M.~Raissi, and F.~Piccialli.
\newblock Scientific {Machine} {Learning} {Through} {Physics}–{Informed}
  {Neural} {Networks}: {Where} we are and {What}’s {Next}.
\newblock {\em J. Sci. Comput.}, 92(3):88, July 2022.

\bibitem{cybenko_approximation_1989}
G.~Cybenko.
\newblock Approximation by superpositions of a sigmoidal function.
\newblock {\em Math. Control Signals Syst.}, 2(4):303--314, December 1989.

\bibitem{dissanayake_neural-network-based_1994}
M.~W. M.~G. Dissanayake and N.~Phan-Thien.
\newblock Neural-network-based approximations for solving partial differential
  equations.
\newblock {\em Commun. Numer. Methods Eng.}, 10(3):195--201, 1994.

\bibitem{dolean2024multilevel}
V.~Dolean, A.~Heinlein, S.~Mishra, and B.~Moseley.
\newblock Multilevel domain decomposition-based architectures for
  physics-informed neural networks.
\newblock {\em Comput. Methods Appl. Mech. Eng.}, 429:117116, 2024.

\bibitem{weinan2017}
W.~E, J.~Han, and A.~Jentzen.
\newblock Deep learning-based numerical methods for high-dimensional parabolic
  partial differential equations and backward stochastic differential
  equations.
\newblock {\em Commun. Math. Stat.}, 5(4):349--380, 2017.

\bibitem{yu2018deep}
W.~E and B.~Yu.
\newblock The {Deep} {Ritz} {Method}: {A} {Deep} {Learning}-{Based} {Numerical}
  {Algorithm} for {Solving} {Variational} {Problems}.
\newblock {\em Commun. Math. Stat.}, 6(1):1--12, March 2018.

\bibitem{heusel2017gans}
M.~Heusel, H.~Ramsauer, T.~Unterthiner, B.~Nessler, and S.~Hochreiter.
\newblock {GANs} trained by a two time-scale update rule converge to a local
  nash equilibrium.
\newblock In {\em Proceedings of the 31st {International} {Conference} on
  {Neural} {Information} {Processing} {Systems}}, {NIPS}'17, pages 6629--6640,
  Red Hook, NY, USA, December 2017. Curran Associates Inc.

\bibitem{howard_stacked_2024}
A.~A. Howard, S.~H. Murphy, S.~E. Ahmed, and P.~Stinis.
\newblock Stacked networks improve physics-informed training: {Applications} to
  neural networks and deep operator networks.
\newblock {\em Found. Data Sci.}, pages 0--0, June 2024.

\bibitem{jacot_neural_2018}
A.~Jacot, F.~Gabriel, and C.~Hongler.
\newblock Neural tangent kernel: convergence and generalization in neural
  networks.
\newblock In {\em Proceedings of the 32nd {International} {Conference} on
  {Neural} {Information} {Processing} {Systems}}, {NIPS}'18, pages 8580--8589,
  Red Hook, NY, USA, December 2018. Curran Associates Inc.

\bibitem{jang2024partitioned}
D.-K. Jang, K.~Kim, and H.~H. Kim.
\newblock Partitioned neural network approximation for partial differential
  equations enhanced with {L}agrange multipliers and localized loss functions.
\newblock {\em Comput. Methods Appl. Mech. Eng.}, 429:117168, 2024.

\bibitem{karniadakis_physics-informed_2021}
G.~E. Karniadakis, I.~G. Kevrekidis, L.~Lu, P.~Perdikaris, S.~Wang, and
  L.~Yang.
\newblock Physics-informed machine learning.
\newblock {\em Nat. Rev. Phys.}, 3(6):422--440, June 2021.

\bibitem{kharazmi_hp-vpinns_2021}
E.~Kharazmi, Zhongqiang Zhang, and George E.~M. Karniadakis.
\newblock $hp$-{VPINNs}: variational physics-informed neural networks with
  domain decomposition.
\newblock {\em Comput. Methods Appl. Mech. Eng.}, 374:Paper No. 113547, 25,
  2021.

\bibitem{adam}
D.~P. Kingma and J.~Ba.
\newblock Adam: A method for stochastic optimization.
\newblock In {\em International Conference on Learning Representations (ICLR)},
  2015.

\bibitem{kissas2020}
G.~Kissas, Y.~Yang, E.~Hwuang, W.~R. Witschey, J.~A. Detre, and P.~Perdikaris.
\newblock Machine learning in cardiovascular flows modeling: predicting
  arterial blood pressure from non-invasive 4{D} flow {MRI} data using
  physics-informed neural networks.
\newblock {\em Comput. Methods Appl. Mech. Engrg.}, 358:112623, 28, 2020.

\bibitem{lagaris_artificial_1998}
I.~E. Lagaris, A.~Likas, and D.~I. Fotiadis.
\newblock Artificial neural networks for solving ordinary and partial
  differential equations.
\newblock {\em IEEE Trans. Neural Networks}, 9(5):987--1000, September 1998.
\newblock Conference Name: IEEE Transactions on Neural Networks.

\bibitem{li2018visualizing}
Hao Li, Zheng Xu, Gavin Taylor, Christoph Studer, and Tom Goldstein.
\newblock Visualizing the loss landscape of neural nets.
\newblock {\em Advances in neural information processing systems}, 31, 2018.

\bibitem{D3M}
K.~Li, K.~Tang, T.~Wu, and Q.~Liao.
\newblock {D3M}: A deep domain decomposition method for partial differential
  equations.
\newblock {\em IEEE Access}, 8:5283--5294, 2020.

\bibitem{lin2020gradient}
T.~Lin, C.~Jin, and M.~Jordan.
\newblock On {Gradient} {Descent} {Ascent} for {Nonconvex}-{Concave} {Minimax}
  {Problems}.
\newblock In {\em Proceedings of the 37th {International} {Conference} on
  {Machine} {Learning}}, pages 6083--6093. PMLR, November 2020.
\newblock ISSN: 2640-3498.

\bibitem{liu1989limited}
D.~C. Liu and J.~Nocedal.
\newblock On the limited memory {BFGS} method for large scale optimization.
\newblock {\em Math. Program.}, 45(1):503--528, 1989.

\bibitem{liu_limited_1989}
Dong~C. Liu and Jorge Nocedal.
\newblock On the limited memory {BFGS} method for large scale optimization.
\newblock {\em Mathematical Programming}, 45(1):503--528, August 1989.

\bibitem{lu_deepxde_2021}
L.~Lu, X.~Meng, Z.~Mao, and G.~E. Karniadakis.
\newblock {DeepXDE}: a deep learning library for solving differential
  equations.
\newblock {\em SIAM Rev.}, 63(1):208--228, 2021.

\bibitem{Lu2021DeepXDE:Equations}
L.~Lu, X.~Meng, Z.~Mao, and G.~E. Karniadakis.
\newblock {DeepXDE: A Deep Learning Library for Solving Differential
  Equations}.
\newblock {\em https://doi.org/10.1137/19M1274067}, 63(1):208--228, 2 2021.

\bibitem{maddu2022inverse}
S.~Maddu, D.~Sturm, C.~L. M{\"u}ller, and I.~F. Sbalzarini.
\newblock Inverse {D}irichlet weighting enables reliable training of physics
  informed neural networks.
\newblock {\em Mach. Learn.: Sci. Technol.}, 3(1):015026, 2022.

\bibitem{mcclenny2023self}
L.~D. McClenny and U.~M. Braga-Neto.
\newblock Self-adaptive physics-informed neural networks.
\newblock {\em J. Comput. Phys.}, 474:111722, 2023.

\bibitem{mishra2022estimates}
S.~Mishra and R.~Molinaro.
\newblock Estimates on the generalization error of physics-informed neural
  networks for approximating a class of inverse problems for {PDEs}.
\newblock {\em IMA J. Numer. Anal.}, 42(2):981--1022, 2022.

\bibitem{muller_achieving_2023}
J.~Müller and M.~Zeinhofer.
\newblock Achieving {High} {Accuracy} with {PINNs} via {Energy} {Natural}
  {Gradient} {Descent}.
\newblock In {\em Proceedings of the 40th {International} {Conference} on
  {Machine} {Learning}}, pages 25471--25485. PMLR, July 2023.
\newblock ISSN: 2640-3498.

\bibitem{Nabian2021EfficientSampling}
M.~A. Nabian, R.~J. Gladstone, and H.~Meidani.
\newblock {Efficient training of physics-informed neural networks via
  importance sampling}.
\newblock {\em Computer-Aided Civil and Infrastructure Engineering},
  36(8):962--977, 4 2021.

\bibitem{nguyen-thanh_deep_2020}
V.~M. Nguyen-Thanh, X.~Zhuang, and T.~Rabczuk.
\newblock A deep energy method for finite deformation hyperelasticity.
\newblock {\em Eur. J. Mech. A. Solids}, 80:103874, March 2020.

\bibitem{paszke_pytorch_2019}
A.~Paszke, S.~Gross, F.~Massa, A.~Lerer, J.~Bradbury, G.~Chanan, T.~Killeen,
  Z.~Lin, N.~Gimelshein, L.~Antiga, A.~Desmaison, A.~Kopf, E.~Yang, Z.~DeVito,
  M.~Raison, A.~Tejani, S.~Chilamkurthy, B.~Steiner, L.~Fang, J.~Bai, and
  S.~Chintala.
\newblock {PyTorch}: {An} {Imperative} {Style}, {High}-{Performance} {Deep}
  {Learning} {Library}.
\newblock In {\em Advances in {Neural} {Information} {Processing} {Systems}},
  volume~32. Curran Associates, Inc., 2019.

\bibitem{rahaman_spectral_2019}
N.~Rahaman, A.~Baratin, D.~Arpit, F.~Draxler, M.~Lin, F.~A. Hamprecht,
  Y.~Bengio, and A.~Courville.
\newblock On the {Spectral} {Bias} of {Neural} {Networks}, May 2019.
\newblock arXiv:1806.08734 [cs, stat].

\bibitem{raissi_physics-informed_2024}
M.~Raissi, P.~Perdikaris, N.~Ahmadi, and G.~E. Karniadakis.
\newblock Physics-{Informed} {Neural} {Networks} and {Extensions}, August 2024.
\newblock arXiv:2408.16806 [cs].

\bibitem{raissi2019}
M.~Raissi, P.~Perdikaris, and G.~E. Karniadakis.
\newblock Physics-informed neural networks: a deep learning framework for
  solving forward and inverse problems involving nonlinear partial differential
  equations.
\newblock {\em J. Comput. Phys.}, 378:686--707, 2019.

\bibitem{shi2021comparative}
E.~Shi and C.~Xu.
\newblock A comparative investigation of neural networks in solving
  differential equations.
\newblock {\em J. Algorithms Comput. Technol.}, 15:1748302621998605, 2021.

\bibitem{sirignano2018}
J.~Sirignano and K.~Spiliopoulos.
\newblock D{GM}: a deep learning algorithm for solving partial differential
  equations.
\newblock {\em J. Comput. Phys.}, 375:1339--1364, 2018.

\bibitem{Son2023}
H.~Son, S.~W. Cho, and H.~J. Hwang.
\newblock Enhanced physics-informed neural networks with augmented {L}agrangian
  relaxation method ({AL-PINN}s).
\newblock {\em Neurocomputing}, page 126424, 2023.

\bibitem{sun2024domain}
Q.~Sun, X.~Xu, and H.~Yi.
\newblock Domain decomposition learning methods for solving elliptic problems.
\newblock {\em SIAM J. Sci. Comput.}, 46(4):A2445--A2474, 2024.

\bibitem{tancik_fourier_2020}
M.~Tancik, P.~P. Srinivasan, B.~Mildenhall, S.~Fridovich-Keil, N.~Raghavan,
  U.~Singhal, R.~Ramamoorthi, J.~T. Barron, and R.~Ng.
\newblock Fourier features let networks learn high frequency functions in low
  dimensional domains.
\newblock In {\em Proceedings of the 34th {International} {Conference} on
  {Neural} {Information} {Processing} {Systems}}, {NIPS} '20, pages 7537--7547,
  Red Hook, NY, USA, December 2020. Curran Associates Inc.

\bibitem{toscano_pinns_2024}
J.~D. Toscano, V.~Oommen, A.~J. Varghese, Z.~Zou, N.~A. Daryakenari, C.~Wu, and
  G.~E. Karniadakis.
\newblock From {PINNs} to {PIKANs}: {Recent} {Advances} in {Physics}-{Informed}
  {Machine} {Learning}, October 2024.
\newblock arXiv:2410.13228.

\bibitem{visser_pacmann_2024}
C.~Visser, A.~Heinlein, and B.~Giovanardi.
\newblock {PACMANN}: {Point} {Adaptive} {Collocation} {Method} for {Artificial}
  {Neural} {Networks}, November 2024.
\newblock arXiv:2411.19632.

\bibitem{wang2023expert}
S.~Wang, Shyam Sankaran, Hanwen Wang, and P.~Perdikaris.
\newblock An expert's guide to training physics-informed neural networks.
\newblock {\em arXiv preprint arXiv:2308.08468}, 2023.

\bibitem{NTK-theory}
S.~Wang, Xinling Yu, and P.~Perdikaris.
\newblock When and why {PINN}s fail to train: {A} neural tangent kernel
  perspective.
\newblock {\em J. Comput. Phys.}, 449:110768, 2022.

\bibitem{wang_multi-stage_2024}
Y.~Wang and C.-Y. Lai.
\newblock Multi-stage neural networks: {Function} approximator of machine
  precision.
\newblock {\em J. Comput. Phys.}, 504:112865, May 2024.

\bibitem{wu2023comprehensive}
C.~Wu, M.~Zhu, Q.~Tan, Y.~Kartha, and L.~Lu.
\newblock A comprehensive study of non-adaptive and residual-based adaptive
  sampling for physics-informed neural networks.
\newblock {\em Comput. Methods Appl. Mech. Eng.}, 403:115671, 2023.

\bibitem{xu_overview_2022}
Z.-Q.~J. Xu, Y.~Zhang, and T.~Luo.
\newblock Overview frequency principle/spectral bias in deep learning, October
  2022.
\newblock arXiv:2201.07395 [cs].

\bibitem{yang_iterative_2024}
H.~J. Yang and H.~H. Kim.
\newblock Iterative algorithms for partitioned neural network approximation to
  partial differential equations.
\newblock {\em Comput. Math. Appl.}, 170:237--259, September 2024.

\bibitem{yang2021b}
L.~Yang, X.~Meng, and G.~E. Karniadakis.
\newblock B-pinns: Bayesian physics-informed neural networks for forward and
  inverse {PDE} problems with noisy data.
\newblock {\em J. Comput. Phys.}, 425:109913, 2021.

\bibitem{Yang2019}
Y.~Yang and P.~Perdikaris.
\newblock Adversarial uncertainty quantification in physics-informed neural
  networks.
\newblock {\em J. Comput. Phys.}, 394:136--152, 2019.

\end{thebibliography}

\end{document}